\documentclass[12pt]{article}

\usepackage{times}
\usepackage{epsfig}
\usepackage{epstopdf}
\usepackage{graphicx}
\usepackage{amsmath}
\usepackage{amssymb}

\usepackage{multirow}
\usepackage{psfrag}

\usepackage{color}
\usepackage{subfigure}

\usepackage{algorithm}
\usepackage{algorithmic}
\usepackage{multirow}

\usepackage[margin=2.5cm]{geometry}

\newcommand{\Real}[0]{\mathbb{R}}

\newcommand{\Complex}[0]{\mathbb{C}}
\newcommand{\Natural}[0]{\mathbb{N}}

\newcommand{\T}[0]{\top}
\newcommand{\gmv}[1]{\mbox{\boldmath$#1$}}
\newcommand{\mv}[1]{\mathbf{#1}}
\newcommand{\mm}[1]{\mathbf{#1}}

\newcommand{\Rep}[0]{\mv D}

\newcommand{\mRep}[1]{\mm{\Rep}_{#1}}

\newcommand{\im}[0]{\mathbf{i}}
\newcommand{\conj}[1]{ \overline{#1} }

\title{Left-Invariant Diffusion on the Motion Group in terms of the Irreducible Representations of $SO(3)$}

\author{Marco~Reisert and Henrik~Skibbe\\\ \\
University Medical Center Freiburg, Medical Physics, Germany\\
\emph{marco.reisert@uniklinik-freiburg.de}
}
\begin{document}


\maketitle

%
%
%
%
%

\begin{abstract}
In this work we study the formulation of convection/diffusion equations on the 3D motion group $SE(3)$ in terms of the irreducible representations of $SO(3)$. 
Therefore, the left-invariant vector-fields on $SE(3)$ are expressed as linear operators, that are differential forms in the translation coordinate and 
algebraic in the rotation. In the context of 3D image processing this approach avoids the explicit discretization of $SO(3)$ or $S_2$, respectively.
This is particular important for $SO(3)$, where a direct discretization is infeasible due to the enormous memory consumption. 
We show two applications of the framework: one in the context of diffusion-weighted magnetic resonance imaging and one in the context of 
object detection. 
\end{abstract}
\newcommand{\Tk}[1]{\mathcal{T}_{#1}}
\newcommand{\Jk}[1]{\mathcal{J}_{#1}}
\newcommand{\Tx}{\mathcal{T}_x}
\newcommand{\Ty}{\mathcal{T}_y}
\newcommand{\Tz}{\mathcal{T}_z}
\newcommand{\Jx}{\mathcal{J}_x}
\newcommand{\Jy}{\mathcal{J}_y}
\newcommand{\Jz}{\mathcal{J}_z}
\newcommand{\Jp}{\mathcal{J}_+}
\newcommand{\Jm}{\mathcal{J}_-}
\newcommand{\Jpm}{\mathcal{J}_{\pm 1}}
\newcommand{\RLP}{\mathcal{J}^2}
\newcommand{\lalg}[0]{se(3)}
\newcommand{\mogrp}[0]{\text{SE}(3)}
\newcommand{\Lsq}[0]{\mathbb L_2}
\newcommand{\lif}[0]{\Xi}
\newcommand{\ME}[0]{Z}

{\footnotesize
\textbf{Keywords:}
Spherical Harmonics, Left-invariant Diffusion, Partial Differential Equations, Wigner D-Matrix, Clebsch-Gordan coefficients, Diffusion weighted MR-imaging, High Angular Resolution Diffusion Imaging (HARDI),
Diffusion Tensor Imaging (DTI), Spatial Regularization, Spherical Hough Transform
}
\section{Introduction}
Image Processing in 3D becomes more and more popular and necessary due to the enormous amount of scientific data acquired with modern imaging 
techniques like magnetic resonance imaging, computer tomography or confocal laser microscopy to mention only a few. Still,
the processing of directional or tensorial information derived from primary modalities or directly measured like in diffusion weighted imaging (DWI) 
drives todays hardware to its limits. For example, imagine a typical DWI measurement with imaging matrix of $100^3$. If we want to sufficiently 
represent the orientation space with e.g. $500$ points, we need already about 4 GB of memory for just one instance. As typical algorithms (e.g. like 
Conjugate Gradients) usually require for more than one instance, we are already at the limits of a common desktop PC. 
In this article we describe how the generators of diffusion and convection on $\Real^3 \times S^2$ and $\Real^3 \times SO(3)$ can
 be described and implemented in terms of the irreducible representations of the 3D rotation group. Most of the implementations 
solving $\Real^3 \times S^2$-diffusion equations \cite{Barm2010,delputte:postprocessing,duits:leftinv,reisert2011} relied on an equiareal discretization of the two-sphere $S_2$,
while implementation for functions $SE(3) \mapsto \Complex$ even do not exist due to the enormous 
memory consumption, although there might be several useful applications like feature detection for non-rotation symmetric templates.
Indeed, there are implementations \cite{duits:leftinv} that use spherical harmonics as an intermediate $S_2$-interpolation
scheme, but they cannot benefit of the well-known advantages of the spherical harmonic representation, like 
the compact and memory efficient storage, the analytic and efficient computations of $S_2$-convolutions,
and the closeness under rotations. The aim of this paper is the formulation of common differential 
operators acting on functions $\Real^3 \times S_2 \mapsto \Complex$ in terms of spherical harmonics.
More generally, we show how a large class partial differential equations on functions $SE(3) \mapsto \Complex$  can be represented
in terms of the irreducible representations of the rotations group $SO(3)$ and solved without any angular discretization. 
In this way one can benefit from all the advantages which the harmonic representations offer. A discretization of the $S_2$ or $SO(3)$, respectively, is avoided, and
one is able to implement diffusion on the full $SE(3)$ with reasonable memory consumption. We show one application in the context 
of high angular resolution diffusion MR-imaging, where spherical harmonic representations are common and compare to the equiareal representation.
Second, it is shown how the proposed framework can be used to implement the spherical Hough transform in an efficient manner. 

\subsection{Related Work}
In the context of line and contour enhancement in 2D the special motion group $SE(2)$ plays a key role \cite{almsick2005,duitsscale,duitsparabolic,weickert:ce}.
It can be used to set up a scale space theory. More recently, extensions to 3D of these concepts appeared \cite{duits:leftinv}. While the applications in 2D are 
typically related to feature detection and image enhancement, the 3D extension offers a new application field: the processing of diffusion weighted magnetic 
resonance images (DWI). In DWI already the acquired measurements are functions on $\Real^3\times S_2$. 
Based on the  directional dependency of water diffusivity in fibrous tissue of the human brain it is possible to reveal underlying connectivity information. 
One of the main challenges in DWI is the estimation in so-called fiber/diffusion orientation distributions. 
There are numerous methods for estimating orientation distributions:
classical Q-ball imaging \cite{Tuch2004}, constrained spherical deconvolution \cite{Tournier2007},
proper probability density estimation \cite{Aganj2010,TVega-NeuroImage09,CanalesMRM2009,Barnett2009}
and  spatially regularized density estimations for tensor-valued images \cite{Goh2009,reisert2011,Tschumperle03,Burgeth2009,Savadjiev2006,Barmpoutis2009}.
Most of the employed algorithms rely on tensorial or spherical harmonic representation of the orientation distributions.
On the other hand, most of the algorithms for orientation distribution estimation that consider the local surrounding of a voxel and use intervoxel information rely
on a discretization of the two-sphere \cite{Barm2010,delputte:postprocessing,duits:leftinv,reisert2011}. 

In two dimensions the representation of orientation and tensor fields in terms of circular harmonics (or, the irreducible representations of $SO(2)$)
is relatively simple and quite frequent in literature \cite{franken:tv,re:tip08,reisert:cdf,claudio:tip}. Complex Calculus offers a well-founded background:
the ordinary Cartesian partial derivatives $\partial_x,\partial_y$ are replaced by the complex ones $\partial_z = (\partial_x - \im \partial_y)/2$ and
$\partial_{\conj{z}} = (\partial_x + \im \partial_y)/2$. In \cite{re:bu:tiipcv09,SphericalDerivativesPAMI} three-dimensional derivative operators are introduced that behave similar to
complex derivatives, that is, they are compliant with the rotation behavior of spherical harmonics in 3D. In \cite{chiri:2000,chiri:book} the Fourier transform 
of $SE(3)$ is used in the context of engineering applications. For the efficient computation of $SE(3)$-convolution functions are expressed in 
terms of the unitary irreducible representations (UiR) of $SE(3)$. The present work proposes a kind of intermediate representation. While 
the full $SE(3)$-UiR representation decomposes also the spatial variable in terms of Bessel-functions, the representation in terms 
of $SO(3)$-UiR leaves the spatial part untouched.

\subsection{Preliminaries and Organization}
Most of the mathematical notations and conventions are adopted from \cite{duits:leftinv} regarding the geometry and parametrization of $SE(3)$ and $SO(3)$. 
Regarding the irreducible representations the notations are similar to \cite{SphericalDerivativesPAMI}. 

In Section \ref{sec:diffgeo} we give a short introduction into the geometry of $SE(3)$ and fix the related conventions. 
Section \ref{sec:irred} introduces the necessary background of representation theory of $SO(3)$. 
The main contribution of this work is presented in Section \ref{sec:main}, where the left-invariant vector-fields of $SE(3)$ 
are expressed in terms of the irreducible representations of $SO(3)$. 
Finally, in Section \ref{sec:exp} applications in the context of DWI and object detection are proposed. 
\section{Life in SE(3)} \label{sec:diffgeo}

The motion group is the semidirect product of the rotation group $SO(3)$ and the translation group. An element $g\in \mogrp$ is composed 
of a translation vector $(x,y,z) = \mv r \in \Real^3$ and a rotation matrix $\mv R\in SO(3)$ with multiplication law
\[
(\mv r, \mv R) (\mv r',\mv R') = (\mv R \mv r' + \mv r, \mv R \mv R')
\]
The Lie group $\mogrp$ is generated by the six dimensional Lie algebra $\lalg = T \mogrp$ spanned by the six
left-invariant vector fields $  \{\Tx,\Ty,\Tz,\Jx,\Jy,\Jz\}$ corresponding to the three translations 
and the three rotation axis. They generate the right regular motion of smooth functions 
$\phi : SE(3) \mapsto \Complex$, i.e. if we define $(\mathcal{R}_h \phi)(g) := \phi(gh)$, then
the application of a left-invariant vector-field $\mathcal{A}$ gives
\[
(\mathcal{A}\phi)(g) = \left. \frac{d}{dt} \right\vert_{t=0}(\mathcal{R}_{h(t)} \phi)(g) = \left. \frac{d}{dt} \right\vert_{t=0} \phi(g h(t))
\]
where $h : \Real \mapsto \mogrp$ is a smooth curve with $h(0)=e$. And hence
\[
\phi(gh) = \mathcal{R}_h\phi(g) = \exp(\mathcal{A}(g) t)\phi(g) 
\]
Note that the tangents $\mathcal{A}(g) \in T_g SE(3)$ do vary over the group and depend on the point $g$. 
%
The vector-fields are usually expressed in Euler angles. 
We parametrized the rotation in Euler-angles $\alpha,\beta,\gamma$ in ZYZ convention as follows
\[
\mv R_g = \mv R_{z,\gamma} \mv R_{y,\beta} \mv R_{z,\alpha}
\]
where all rotations are counter-clockwise. In this parametrization the left-invariant vector-fields
of the translation take the form:
\[
\left( \begin{array}{c} \Tx \\ \Ty \\ \Tz \end{array} \right) =  \mv R_g^\T \nabla 
\]
with $\nabla = (\partial_x,\partial_y,\partial_z)^\T$
and for the rotation
\begin{eqnarray*}
\Jx &=& \cos \alpha \cot \beta \partial_\alpha + \sin \alpha \partial_\beta - \frac{\cos \alpha}{\sin \beta} \partial_\gamma \\
\Jy &=& - \sin \alpha \cot \beta \partial_\alpha + \cos \alpha \partial_\beta + \frac{\sin \alpha}{\sin \beta} \partial_\gamma \\
\Jz &=& \partial_\alpha
\end{eqnarray*}
Note that the fields depend explicitly on the position in the group. The conventions are the same to ones used in \cite{duits:leftinv}.

\section{Unitary Irreducible Representations of SO(3)} \label{sec:irred}
It is well known that the spherical harmonics (see Appendix \ref{ap:sh} for definition) are a orthogonal basis for the square integrable 
functions on the unit sphere $\Lsq(S_2)$. They have a variety of gentle properties with respect
to rotations. The Wigner D-matrices play the same role for square integrable functions $\Lsq(SO(3))$
on the rotation group $SO(3)$ itself. They are representations of $SO(3)$ and are written in
Euler angles as
\begin{eqnarray}
(\mv D^j(g))_{n m} = D^j_{n m}(\gamma,\beta,\alpha) = e^{-\im n \gamma } d^j_{n m} (\beta) e^{-\im m \alpha} \label{eq:dwigexplicit}
\end{eqnarray}
where the real 'small' d-matrix $d^j_{n m}$ is related to the Jacobi polynomials (see Appendix \ref{ap:wigner}).
As the $\mv D^j(g)$ are representations of $SO(3)$ they obey the multiplication law
\[
\mv D^j(g h) = \mv D^j(g) \mv D^j (h)
\]
for any $g,h \in SO(3)$. For each order $j\in \Natural$ they work on complex $2j+1$-dimensional vector spaces, i.e. $-j \leq n,m \leq j$.
  Furthermore, they are the irreducible representations of $SO(3)$, i.e. there is no linear transformation $\mv A$ such that 
$\mv A^\T \mv D^j(g) \mv A$ is block-diagonal for all $g\in SO(3)$. 
For $j=1$ there is the direct relation to the original rotation matrix $\mv R_g$ by
\[
\mv D^1(g) = \mv S \mv R_g \mv S^\T \text{ with } \mv S = \frac{1}{\sqrt{2}} \left( \begin{array}{ccc} -1 & \im & 0 \\ 0 & 0 & \sqrt{2} \\ -1 & -\im &0 \end{array} \right)
\]
Note that $\mv S$ is unitary, i.e. $\mv S^\T \mv S = \mv I$. 
The irreducibility has an important consequence. By the Peter-Weyl theorem the irreducible representations form
a complete orthogonal basis set of functions with respect to the group dot-product:
\[
\langle  D^j_{nm}, D^{j'}_{n'm'} \rangle_{SO(3)} = \int_{g\in SO(3)} D^j_{nm}(g) \conj{D^{j'}_{n'm'}}(g) = \frac{8 \pi^2}{2 j + 1} \delta_{m'm} \delta_{n'n} \delta_{j'j}
\]
Thus, we can write any $\phi \in \Lsq(SO(3))$ as
 \[
\phi(g) = \frac{1}{8\pi^2} \sum_{j=0}^\infty \sum_{n=-j}^j \sum_{m=-j}^j (2j+1)\  \conj{D^j_{nm}}(g)\ f^j_{nm}
\]
where the expansion coefficients can be obtained by a simple projection
\[
 f^j_{nm} = \langle  \phi, \conj{D}^{j}_{nm} \rangle_{SO(3)} = \int_{g\in SO(3)}  \phi(g)\  D^j_{nm}(g)
\]
onto the Wigner D-matrix.
\subsection{Clebsch Gordan Coefficients of SO(3)}
The Clebsch Gordan (CG) coefficients interrelate the irreducible representations of different order. 
We denote a Clebsch Gordan coefficient by $\langle l m | l_1 m_1, l_2 m_2 \rangle$ where $l,l_1,l_2$ 
are different orders such that the triangle inequality $|l_1-l_2| \leq l \leq |l_1+l_2|$ holds, otherwise the coefficients vanish. 
The basic equation connecting two different representations of order $l_1$ and $l_2$ is the following 
\begin{eqnarray}
D^\ell_{mn} = \sum_{m_1+m_2 = m \atop n_1+n_2 = n} D^{\ell_1}_{m_1 n_1}D^{\ell_2}_{m_2 n_2}
 \langle l m | l_1 m_1, l_2 m_2 \rangle
 \langle l n | l_1 n_1, l_2 n_2 \rangle \label{dw:decomp1}
\end{eqnarray}
Note the additional selection rule of the CG-coefficients, they only contribute if $m = m_1+ m_2$. One also knows
that  $\langle j 0 | j_1 0, j_2 0 \rangle = 0$ if $j+j_1+j_2$ is odd. 
There is a variety of orthogonality and symmetry relations for the Clebsch-Gordan coefficients making them itself 
a similarity transformation. We mention here the most important orthogonality relations:
\begin{eqnarray}
\sum_{j=0}^\infty \sum_{m=-j}^j \langle j m | j_1 m_1, j_2 m_2 \rangle \langle j m | j_1 m'_1, j_2 m'_2 \rangle
&=& \delta_{m_1,m_1'} \delta_{m_2,m_2'}  \label{cg:orth1} \\
\sum_{{m = m_1 + m_2} } \langle j m | j_1 m_1, j_2 m_2 \rangle \langle j' m' | j_1 m_1, j_2 m_2 \rangle
&=& \delta_{j,j'} \delta_{m,m'}\label{cg:orth2}  
\end{eqnarray}
For example, suppose a sequence of numbers $a^j_m$ with $ j \leq J$ is given. Then, for two fixed $j_1,j_2$ 
obeying  $|j_1-j_2| \leq j \leq |j_1+j_2|$ we can compute $f^{j_1,j_2}_{m_1,m_2} = \sum_{j,m} a^j_m \langle j m | j_1 m_1, j_2 m_2 \rangle$, which
contains all the information about the original $a^j_m$ in an unitary way. 
Besides the orthogonality relations there are numerous other symmetry and associativity relations for the Clebsch Gordan 
coefficients, some of them are listed in Appendix \ref{ap:cg}.

\subsection{Solid Harmonics and Spherical Derivatives}
In terms of the associated Legendre polynomials the components $(\mv Y^j)_m = Y^\ell_m$ of the Racah-normalized spherical harmonics  (for further details see Appendix \ref{ap:sh}) are written as
\[
Y^j_m (\beta,\gamma) = \sqrt{\frac{(j-m)!}{(j+m)!}} P^j_m(\cos(\beta)) e^{\im m\gamma}
\]
Instead of $\beta,\gamma$ defining a point on the sphere, we write in the following $\mv n \in S_2$ as normalized Cartesian vector. 
There is a close relation between the Wigner D-matrix and the spherical harmonics:
\[
\mv D^j(g) \mv Y^j(\mv R_g^\T\mv n) = \mv Y^j(\mv n)
\]
for any $g\in SO(3)$ and $\mv n \in S_2$. That is, the expansion coefficients of a spherical harmonic expansion
rotate by the application of Wigner D-matrices. 
With this we can identify the central column of a Wigner D-matrix with a conjugate spherical harmonic by rotating 
a spherical harmonic along the z-axis
\[
\mv D^j(g) \mv Y^j(\mv e_z) = \mv Y^j(\mv R_g \mv e_z)
\]
and with the additional knowledge that $(\mv Y^j(\mv e_z))_m = \delta_{m,0}$ we have 
\[
D^j_{n0}(\gamma,\beta,0) = \conj{Y^j_n}(\beta,\gamma).
\]
Next to the spherical harmonics, the so called solid harmonics 
\[
R^j_m(\mv r) = |\mv r|^j \mv Y^j_m(\mv r / |\mv r|)
\]
are solutions of the homogeneous Laplace equation, i.e. $\Delta R^j_m = 0$. They are homogeneous polynomials of degree $j$, that is, $R^j_m(\lambda \mv r) = \lambda^j R^j_m(\mv r)$ 
for any $\lambda \in \Real$. So, we can define a differential operator as follows
\[
\gmv \partial^j_m := R^j_m(\nabla)
\]
which is a spherical tensor operator, i.e. it inherits all rotation properties of irreducible representations.
Further note that $\gmv \partial^1 = \mv S \nabla$.

\section{Diffusion Equations in terms of the Irreducible Representations of SO(3)} \label{sec:main}
The quadratic forms in the left-invariant vector fields  \cite{duits:leftinv}  generate the 'dynamics'
on functions $\phi : SE(3) \mapsto \Complex$. Our goal is to find them in terms of 
the irreducible representations of $SO(3)$. More precisely, suppose we have an evolution equation
\[
\partial_t \phi(g,t) = H \phi(g,t)
\]
where $g$ is parametrized as proposed above and the evolution generator is a polynomial in the left-invariant vector fields $H = H(\vec{\mathcal T},\vec{\mathcal J}) = H(\Tx,\Ty,\Tz,\Jx,\Jy,\Jz)$.
That is, $H$ acts as  a differential operator in $\partial_x,\partial_y,\partial_z$ and $\partial_\alpha,\partial_\beta,\partial_\gamma$
on the function $\phi$. By decomposing $\phi$ in terms of Wigner D-matrices
\[
\phi(g) = \phi(g_{\mv r} g_R) = \frac{1}{8\pi^2} \sum_{j=0}^\infty \sum_{n=-j}^j \sum_{m=-j}^j (2j+1)\  \conj{D^j_{nm}}(g_R)\ f^j_{nm}(\mv r,t)
\]
we will be able to show that the evolution equation in terms of $f^j_{nm}(\mv r,t)$ is just 
a differential operator in the spatial coordinates and algebraic for the angular coordinates. The form of the equation will not depend on the
particular choice of the chart chosen to parametrize $SO(3)$. That is, we expect the equation to be
\[
\partial_t f^j_{nm}(\mv r,t) = \sum_{j',n',m'}  \hat{H}^{jnm}_{j'n'm'}\ f^{j'}_{n'm'}(\mv r,t)
\]
where the $\hat{H}^{jnm}_{j'n'm'}$ are differential operators in the spatial coordinates.
To find $\hat{H}$ we have to express all appearing quantities in terms of the irreducible representations $D^j_{nm}$, i.e. 
we have to evaluate the matrix elements $\langle H \conj{D}^j_{nm}, \conj{D}^{j'}_{n'm'} \rangle_{SO(3)} := \hat{H}^{jnm}_{j'n'm'}$ or 
 $\langle  H \phi , \conj{D}^j_{nm}\rangle_{SO(3)}$, respectively, which is subject of the next section. Note, that by the help of 
the UiR of the full group $SE(3)$ \cite{chiri:2000}, the above equation can be made purly algebraic as long as no external gauge field is used. 
However, this approach would also approximate the translation/spatial part of the function, which is in our context not neccessary and
avoids any approximation artefacts coming from this side. 
Of course, the Wigner-D expansion is restricted to a finite cutoff index $j\leq J$ in practice, and thus, we also have approximation errors.
More precisely, if $P_J$ is the linear orthogonal projector onto the subspace 
of functions spanned by the Wigner-D matrices up to an order of $j\leq J$, then the approximated evolution generator is 
$H_J = P_J H P_J$ and the corresponding evolution operator is $\exp(H_Jt) = \exp(P_J H P_Jt)$. That is, the obtained time evolution 
is not the projection $P_J\exp(Ht)P_J$ of the original one, which is important to note.

\subsection{The Left-invariant Vector Fields}
We start with switching to the complex representation which
is obstructed by the irreducible representation of rank $j=1$: 
\begin{eqnarray}
\vec {\mathcal T} = \left( \begin{array}{c} \Tk{-1} \\ \Tk{0} \\ \Tk{+1} \end{array} \right) &=&
\left( \begin{array}{c} -(\Tx - \im\Ty)/\sqrt 2\\ \Tz \\ -(\Tx+\im \Ty)/\sqrt 2 \end{array} \right) = \mv S \mv R_g^\T \nabla =  \mv D^1(g)^\T \gmv \partial^1 \label{eq:transfield1}
\end{eqnarray}
\begin{eqnarray}
\vec{\mathcal J} = \left( \begin{array}{c} \Jk{-1} \\ \Jk{0} \\ \Jk{+1} \end{array} \right) &=&
\left( \begin{array}{c} -(\Jx - \im\Jy)/\sqrt 2\\ \Jz \\ -(\Jx+\im \Jy)/\sqrt 2 \end{array} \right) \label{eq:rotfield1},
\end{eqnarray}
where we formally discriminate between both representations by numeric or character $x,y,z$ subindices. 
The $\Jk{\pm 1}$ are well known from quantum mechanics and engineering literature \cite{chiri:book, wormer:AngMom,miller:harmana,edmonds,rose:angmom,tinkham:grouptheory}
as the body-fixed rigid rotor angular momentum operators. They obey the following equations: 
\begin{eqnarray}
\Jz D^j_{nm} &=& -\im m D^j_{nm} \label{eq:DwigQM1} \\
\Jk{\pm 1} D^j_{nm} &=& \im \sqrt{j(j+1)/2-m(m\pm 1)/2} D^j_{n(m\pm 1)} \label{eq:DwigQM2}.
\end{eqnarray}
Thus, we can directly compute the action of the operators $\Jk{}$ onto the coefficients fields $f^j_{nm}(\mv r)$. 
Therefore, we denote the collection of all coefficients by a bold Latin letter $\mv f$ and we access elements of $\mv f$ by round brackets, i.e. 
$(\mv f)^j_{nm} = f^j_{nm}$ and thus, the coefficient $(\Jz\mv f)^j_{nm}$ is the projection of $(\Jz \phi)$ onto the orthogonal Wigner D-matrix basis:
\[
(\Jz \mv f)^j_{nm} = \langle\Jz \phi,  \conj{D}^j_{nm} \rangle_{SO(3)} = \int_{SO(3)} dg_R\ (\Jz \phi)(g_{\mv r} g_R)\ D^j_{nm}(g_R).
\]
where the integration ranges over $SO(3)$. With formula \eqref{eq:DwigQM1} and partial integration we can find
\begin{eqnarray}
(\Jz \mv f)^j_{nm} &=& -\int dg_R\ \phi(g)\ \Jz D^j_{nm}(g_R)  = \im m\ f^j_{nm} \label{eq:Jrotz}
\end{eqnarray}
and similarly with  equation \eqref{eq:DwigQM2}  we can proceed as follows:
\begin{eqnarray}
(\Jpm \mv f)^j_{nm} &=& -\int dg_R\ \phi(g)\ \Jpm D^j_{nm}(g_R) \nonumber \\
&=& -\im \sqrt{j(j+1)/2-m(m\pm 1)/2} \int dg_R\ \phi(g)\  D^j_{n(m\pm 1)}(g_R) \nonumber \\
&=& -\im \sqrt{j(j+1)/2-m(m\pm 1)/2} \  f^j_{n(m\pm 1)}  \label{eq:Jrotpm}
\end{eqnarray}
The translations are a bit more intricate to compute. We apply the right hand side of equation \eqref{eq:transfield1}, i.e.
$\vec {\mathcal T} = \mv D^1(g)^\T \gmv \partial^1$ onto the field $\phi$ and project onto the orthogonal Wigner basis:
\begin{eqnarray*}
(\Tk{k} \mv f)^j_{nm} &=& \int dg_R\ (\Tk{k}\phi)(g)\  D^j_{nm}(g_R)\\
 &=& \sum_{q=-1,0,1}\int dg_R\  \conj{D^1_{q k}} \gmv \partial^1_{q} \phi(g) D^j_{nm}(g_R).
\end{eqnarray*}
By using the integral formula \eqref{eq:tripleproducts} for triple products of Wigner D-matrices we get
\begin{eqnarray}
(\Tk{k} \mv f)^j_{nm} &=&  \frac{1}{8\pi^2} \sum_{j',n',m',q} (2j'+1) \gmv \partial^1_{q} f^{j'}_{n'm'}  \int dg_R\  \conj{D^1_{q k}} \  \conj{D^{j'}_{n'm'}}\  D^j_{nm} \nonumber \\
 &=&  \sum_{j'=j-1,j,j+1\atop q  = -1,0,1} \sum_{n= n'+q \atop m = m' + k}  \frac{2j'+1}{2 j + 1}  \langle j n | j' n', 1 q \rangle \langle j m | j' m' ,1 k \rangle\ \gmv \partial^1_{q} f^{j'}_{n'm'}\label{eq:transfield2}
\end{eqnarray}
which is already our first main result. It gives the action of the translational left-invariant vector-fields $\Tk{k}$ in terms of the Wigner expansion coefficients $f^j_{nm}$.
Note that the sums are very sparse, which is from a computational viewpoint quite important. The sum for $j'$ runs only over three terms.
\subsection{Quadratic Forms}
While the linear forms in the left-invariant vector-fields generate Euclidean motion, the quadratic forms generate diffusion. 
The important Laplace Beltrami Operator on SO(3)
\begin{equation}
\RLP = \Jx^2 + \Jy^2 + \Jz^2 = \Jk{+1}\Jk{-1} + \Jk{0}^2 + \Jk{-1}\Jk{+1}
\end{equation}
generates diffusion on $SO(3)$. The action of $\RLP$ onto the coefficients $f^j_{nm}$ can be computed with the 
help of equation  \eqref{eq:Jrotz} and \eqref{eq:Jrotpm} to
\begin{equation}
(\RLP \mv f)^j_{n m} = - j(j+1) f^j_{n m}.
\end{equation}
Other arbitrary products $\conj{\Jk{k'}}\Jk{k}$ can be computed quite easily by the use of equation  \eqref{eq:Jrotz} and \eqref{eq:Jrotpm}.
Also products of rotations $\Jk{k}$ and translations $\Tk{k'}$ can be computed by the use of the first-order equations  \eqref{eq:DwigQM1}, \eqref{eq:DwigQM2}, \eqref{eq:transfield2}. 
You can find the detailed formulas in  Appendix \ref{ap:mixedquad}.
However, the product of two translations is again  more cumbersome to evaluate. The vector field $\conj{\Tk{k'}} \Tk{k}$ transforms with respect to $k$ and $k'$ 
like a Cartesian rank 2 tensor. By using equation \eqref{eq:shprod4}  we can write it in terms of spherical tensors as follows (for a proof see Appendix \ref{ap:proofshquad}):
\begin{equation} 
\conj{\Tk{k'}} \Tk{k} =  \frac{\Delta}{3}  -  \frac{\sqrt{10}}{3} \sum_{p=-2}^2  \langle 1 k' | 2 p, 1 k\rangle  \ ((\mv D^2)^\T \gmv \partial^2)_p\label{eq:shquadrel}
\end{equation}
where $\Delta = \partial_x^2 + \partial_y^2 + \partial_z^2$, which is related to the trace of $\conj{\Tk{k'}} \Tk{k}$, i.e. $\Delta = \sum_k |\Tk{k} |^2 $. 
The second term is related to traceless matrix $\conj{\Tk{k'}} \Tk{k}  - \Delta/3 \delta_{kk'}$. 
By using equation \eqref{eq:shquadrel} we can proceed like in the linear case:
\begin{eqnarray*}
(\conj{\Tk{k'} }\Tk{k} \mv f)^j_{nm} &=& \int dg_R\ (\conj{\Tk{k'} }\Tk{k} \phi)(g)\  D^j_{nm}(g_R)\\
&=&  \frac{\Delta}{3} f^j_{nm}  -  \frac{\sqrt{10}}{3\cdot 8 \pi^2}   \sum_{q,p=-2,\hdots,2 \atop j',n',m'} (2j'+1) \gmv \partial^2_q  f^{j'}_{n'm'}    \langle 1 k' | 2 p, 1 k\rangle  \int  \conj{D^2_{q p}} 
\  \conj{D^{j'}_{n'm'}} 
D^j_{nm}\\
&=&  \frac{\Delta}{3}f^j_{nm}  -   \frac{\sqrt{10}}{3}    \sum_{{k' = p + k} \atop {{j'=j-2,\hdots,j+2} \atop {n=n'+q,m=m'+p}}} \frac{2j'+1}{2j+1}    \langle 1 k' | 2 p, 1 k\rangle 
\langle j n | j' n', 2 q \rangle \langle j m | j' m', 2 p\rangle  \gmv \partial^2_q  f^{j'}_{n'm'}
\end{eqnarray*}
where we have again used formula \eqref{eq:tripleproducts} for the triple products of Wigner D-matrices. 
To simplify the formula we consider the special case $k=k'=0$, i.e. $|\Tk{0}|^2 = \Tz^2$, which gives
\begin{eqnarray}
(\Tz^2 \mv f)^j_{nm} &=&  \frac{\Delta}{3} f^j_{nm}  +  \frac{2}{3}   \sum_{ {{j'=j-2,\hdots,j+2} \atop {n=n'+q}}} \frac{2j'+1}{2j+1} 
\langle j n | j' n', 2 q \rangle \langle j m | j' m, 2 0\rangle  \gmv \partial^2_q  f^{j'}_{n'm}.
\end{eqnarray}
On the other hand  $|\Tk{-1}|^2 + |\Tk{1}|^2 =  \Tx^2 + \Ty^2$ gives
\begin{eqnarray}
((\Tx^2 + \Ty^2)\mv f)^j_{nm} &=& \frac{2\Delta}{3} f^j_{nm}  -  \frac{2}{3}   \sum_{ {{j'=j-2,\hdots,j+2} \atop {n=n'+q}}} \frac{2j'+1}{2j+1} 
\langle j n | j' n', 2 q \rangle \langle j m | j' m, 2 0\rangle  \gmv \partial^2_q  f^{j'}_{n'm}.
\end{eqnarray}
Note that there is only a small change between $\Tx^2 + \Ty^2$ and $\Tz^2$. One can easily see that 
$\Tz^2 + \Tx^2 + \Ty^2 = \Delta$. 

\subsection{Restriction to  $\Real^3 \times S_2$}
\begin{figure}
\centering
\includegraphics[width=14cm]{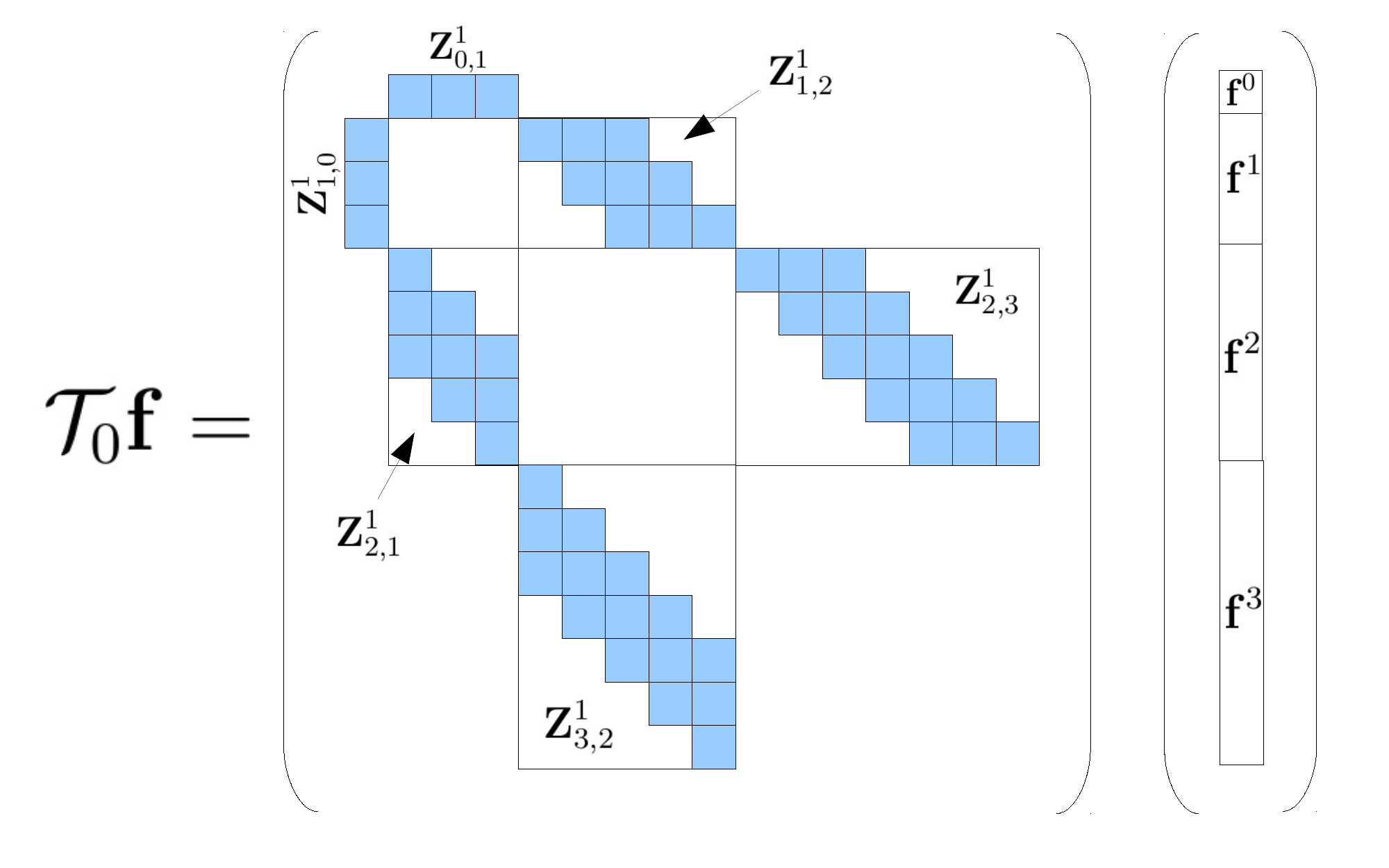} \label{fig:transmat}
\end{figure}

Given a function $\phi : \Real^3 \times S_2 \mapsto \Complex$  it is natural to extend to functions $\hat \phi : SE(3) \mapsto \Complex$ by
\[
\hat \phi(\mv r, \mv R) = \phi(\mv r, \mv R \mv e_z) = \phi(\mv r, \mv n) 
\]
Hence, there is the symmetry $\hat \phi(\mv r, \mv R_{z,\alpha}) = \hat \phi(\mv r, \mv I)$ for any rotation around the z-axis. In Euler angles 
this means that $\hat \phi(\mv r, \gamma, \beta, \alpha) =\hat \phi(\mv r, \gamma, \beta, 0)$. 
By this restriction the evolution operator $H$ has to be invariant with respect to a rotation around the z-axis
\[
H(\vec {\mathcal T},\vec {\mathcal J}) = H(\mv R_{z,\alpha}\vec {\mathcal T},\mv R_{z,\alpha} \vec {\mathcal J}) 
\]
(see \cite{duits:leftinv} for details). Note that such an invariant $H$ always depends on $\Tz,\Tx^2 + \Ty^2$ and $\Jx^2 + \Jy^2$ (at least if no other 
external gauge image is used). To understand the implications for the Wigner D-expansion we have to look for 
Wigner D-matrices that satisfy $D^j_{nm}(\gamma,\beta,0) = D^j_{nm}(\gamma,\beta,\alpha)$ for any $\alpha$. From the explicit form in equation \eqref{eq:dwigexplicit}
it is easy to see that this holds for all $D^j_{n0}$. So, the expansion in terms of Wigner D-matrices reduces to the usual spherical harmonic expansion
\begin{eqnarray*}
\phi(\mv r,\gamma,\beta) &=&  \hat \phi(\mv r, \gamma,\beta,0) = \frac{1}{8\pi^2} \sum_{j=0}^\infty \sum_{n=-j}^j \sum_{m=-j}^j (2j+1)\  \conj{D^j_{nm}}(\gamma,\beta,0)\ f^j_{n0} \\
&=& \frac{1}{8\pi^2} \sum_{j=0}^\infty \sum_{n=-j}^j (2j+1)\  Y^j_{n}(\gamma,\beta)\ f^j_{n}
\end{eqnarray*}
To translate the general action of the translational vector-fields in equation \eqref{eq:transfield2} onto this special case, we just have to set $m=m'=0$. 
The selection rules of the Clebsch-Gordan coefficients imply $k=0$, meaning that only the $\Tk{0}$ gives non-vanishing contributions, which 
is in agreement with \cite{duits:leftinv} that only those operators survive the construction that are well defined on the cosets $SO(3)/SO(2)$. 
For this special case we visualize the matrix elements graphically in Figure \ref{fig:transmat}, where we use the abbreviation
\[
\ME^{J}_{jj',n n'} = (\mv \ME^{J}_{jj'})_{n n'} = \sum_{q  = -J,\hdots,J} \delta_{n,n'+q} \frac{2j'+1}{2 j + 1}  \langle j n | j' n', J q \rangle \langle j 0 | j' 0 ,J 0 \rangle\ \gmv \partial^
J_{q} 
\]
for the block matrices on the upper and lower secondary block diagonal, where $J$ is the differential order of the operator. The rank $j$ is the rank of spherical tensor field
which is returned after application of $\mv \ME^{J}_{jj'}$ and $j'$ is the rank of the incoming spherical tensor field. 
So, we can write equation \eqref{eq:transfield2} for $k=0$ more compactly as
\begin{eqnarray}
(\Tk{0} \mv f)^j_{n}  &=&  \sum_{n'=-(j+1)}^{j+1} \ME^1_{j(j+1),n n'}\ f^{j+1}_{n'} +  \sum_{n'=-(j-1)}^{j-1}  \ME^1_{j(j-1),n n'}\ f^{j-1}_{n'} \nonumber \\
  &=&  \left(\mv \ME^1_{j(j+1)} \mv f^{j+1} + \mv \ME^1_{j(j-1)} \mv f^{j-1} \right)_n \label{eq:transinZ}
\end{eqnarray}
Note that the matrices $\mv \ME$ are also very sparse, the non-zero elements are shaded in blue in Figure \ref{fig:transmat}.
In the same fashion we can write the second order operators 
\begin{eqnarray}
(\Tz^2 \mv f)^j &=&  \frac{\Delta}{3} \mv f^j  +  \frac{2}{3}   ( \mv \ME^2_{j(j+2)} \mv f^{j+2} + \mv \ME^2_{jj} \mv f^j+ \mv \ME^2_{j(j-2)} \mv f^{j-2}) \label{eq:T2res}
\end{eqnarray}
and 
\begin{eqnarray}
((\Tx^2 + \Ty^2) \mv f)^j &=&  \frac{2\Delta}{3} \mv f^j  -  \frac{2}{3}   ( \mv \ME^2_{j(j+2)} \mv f^{j+2}+ \mv \ME^2_{jj} \mv f^{j}+ \mv \ME^2_{j(j-2)} \mv f^{j-2}) 
\end{eqnarray}
Both operators generate horizontal and orthogonal diffusion (\cite{duits:leftinv}). There are also generalizations of these diffusion generators by convolutions
with rotation symmetric functions on $S_2$. These convolutions are known to have a very simple diagonal form in terms of spherical harmonics. If we denote the convolution
operator by $\mv C$, then the application is defined to be $(\mv C \mv f)^j = c_j \mv f^j$, where $c_j$ are the expansion coefficients of the rotation symmetric convolution kernel on $S_2$. 
Two operators, which can be derived from a Tikhononv-regularization problem are $\mv C^\T \Tz^2 \mv C$ and on the other hand $\Tz \mv C^\T \mv C \Tz$. They can easily expressed 
in our framework:
\begin{eqnarray}
(\mv C^\T \Tz^2 \mv C\ \mv f)^j =  a_1 \frac{\Delta}{3} \mv f^j   +   \frac{2}{3}   (a_2 \mv \ME^2_{j(j+2)} \mv f^{j+2}  + a_3\mv \ME^2_{jj} \mv f^j + a_4 \mv \ME^2_{j(j-2)} \mv f^{j-2} ) \label{eq:T2resC}
\end{eqnarray}
where the factors are  $a_1 = a_3 = |c_j|^2$,$a_2 = \conj{c_j}c_{j+2}$ and $a_4 = \conj{c_j} c_{j-2}$. On the other we have 
\begin{eqnarray}
(\Tz \mv C^\T \mv C\Tz  \mv f)^j =  a_1 \frac{2\Delta}{3} \mv f^j   +  \frac{2}{3}   (  a_2 \mv \ME^2_{j(j+2)} \mv f^{j+2}  + a_3 \mv \ME^2_{jj} \mv f^j  +  a_4 \mv \ME^2_{j(j-2)} \mv f^{j-2} ) 
\end{eqnarray}
where $a_1 = a_3 = (|c_{j+1}|^2 +|c_{j-1}|^2)$, $a_2 = |c_{j+1}|^2$ and $a_4 = |c_{j-1}|^2$. Note that the same operators can also be derived for
orthogonal diffusion $(\Tx^2 + \Ty^2)$.

\section{Applications and Experiments} \label{sec:exp}
We want to show two examples where the introduced framework can be applied. First, an example in MR-imaging will show that the introduced operators 
can effectively used as regularization terms for the estimation of fiber orientation distributions on the basis of HARDI-imaging. In this experiment the
SH-based approach will be compared with its discrete counterpart. 
And secondly, we use the translation operator to implement a 3D circular 
Hough transform efficiently. But before starting, some details on the discrete implementation of the derivative operators are given and some general issues are discussed.

\subsection{Discrete Implementation}
All differential operators are implemented by ordinary finite difference schemes. To keep the implementation efficient we restricted our implementation to 
first order approximations.
For the discretization of the first-order derivative operators $\mv \ME^1_{jj'}$ we use a simple central finite difference approximation, that is, the convolution kernel of 
a partial derivative along an arbitrary axis reads $[-1\ 0\ 1]/2$. We also explicitly implemented second order derivatives $\mv \ME^2_{jj'}$ by the common difference scheme:
For example, the convolution kernels of $\partial_{xx}$ and $\partial_{xy}$ are 
\begin{equation}
d_{xx} = \left[ \begin{array}{ccc} 0 &0 &0 \\ 1 & -2 & 1 \\ 0 & 0 & 0 \end{array} \right], \hspace{0.2cm}d_{xy} = \frac{1}{4}\left[ \begin{array}{ccc} 1 &0 &-1 \\ 0 & 0 & 0 \\ -1 & 0 & 1 \end{array} \right],
\end{equation}
the others can be obtained by permutation and symmetry. From these rather rough approximations we cannot expect too much. In particular, the forward Euler-integration of first-order central 
differences is usually a no-go. Unfortunately, usual solution for convection dominated problems by upwind/downwind schemes cannot be applied due to the global nature of 
the harmonic representations. Nevertheless, we want to rely on these crude approximations for the sake practicability.
Of course, there are several more sophisticated approximation schemes (\cite{weickert:ivcip,kroon2010,duits:diffconv20011}), but due to the inherent slow-down of the computation speed, in particular in 3D, we only 
consider the most simple scheme.  And we will see in the experiments that it is enough to obtain reasonable good results.

\begin{figure}
\centering
\includegraphics[width=8cm]{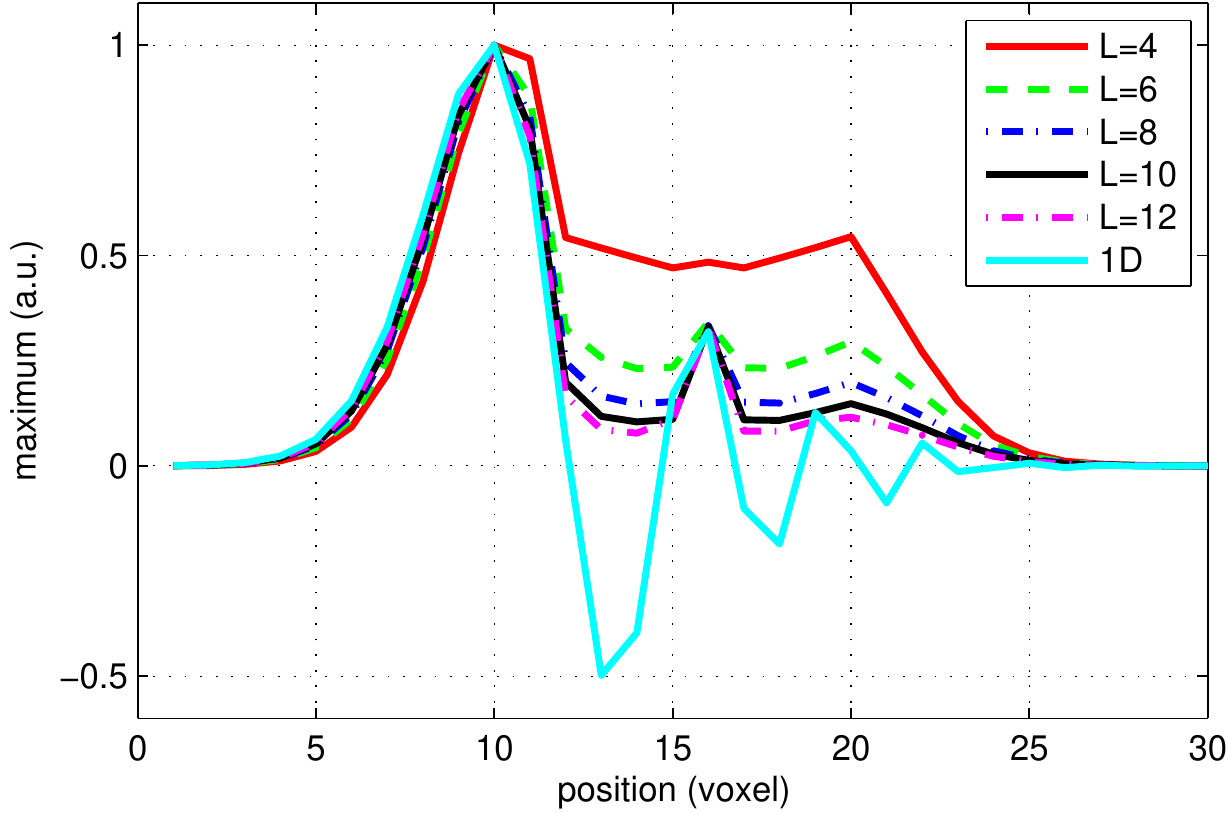}
\includegraphics[width=8cm]{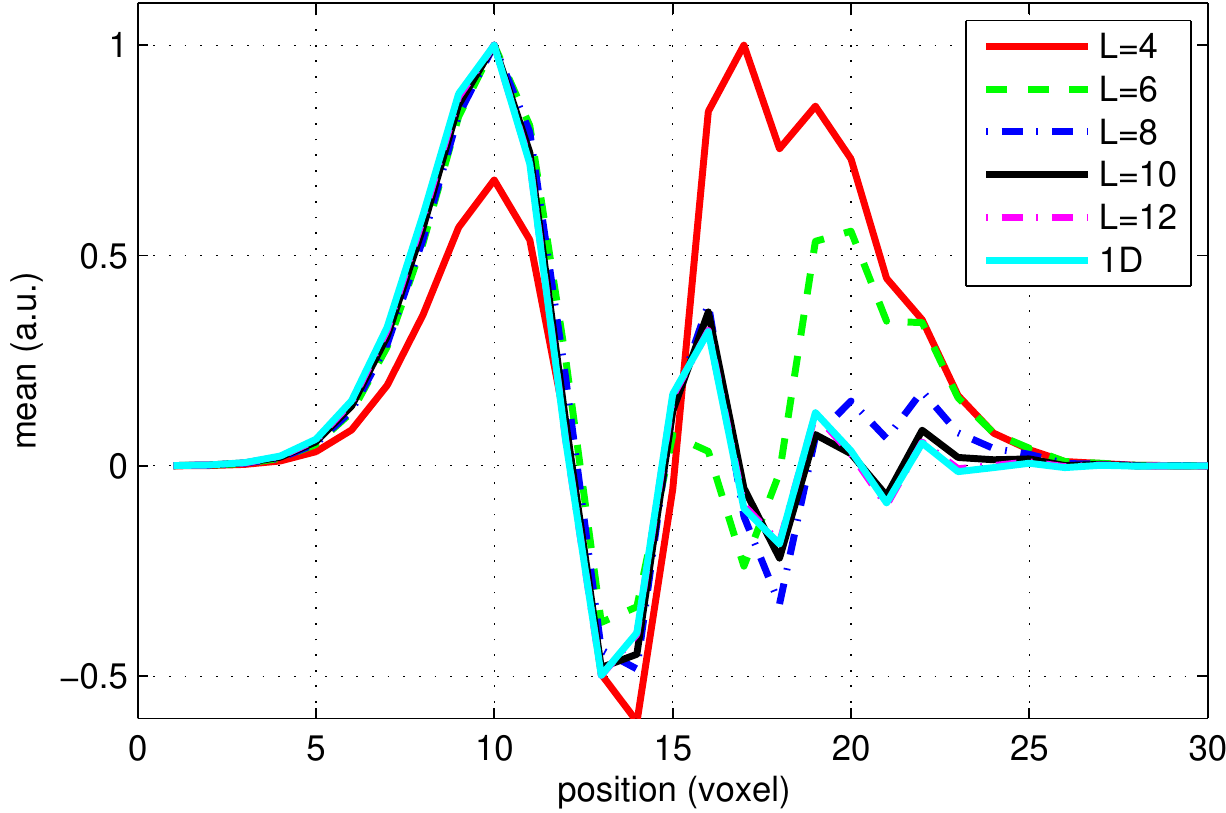} 
\caption{The kernel $\Tk{0}$ was iterated $150$ times with a step-width of $0.05$ resulting in a translation of $7.5$ voxels. The initial condition was 
set to $\phi_0(\mv r,\mv n) = e^{-|\mv r|^2/2} \delta_{\mv n_z}(\mv n)$. 
On the left the maximum  $\max_{\mv n} \phi(z \mv n_z, \mv n)$ along the z-axis is plotted. On the right the mean-component $j=0$, i.e. $\mv f^0(z\mv n_z)$, 
is visualized. } \label{fig:translation_test1}
\end{figure}

\begin{figure}
\centering
\includegraphics[width=8cm]{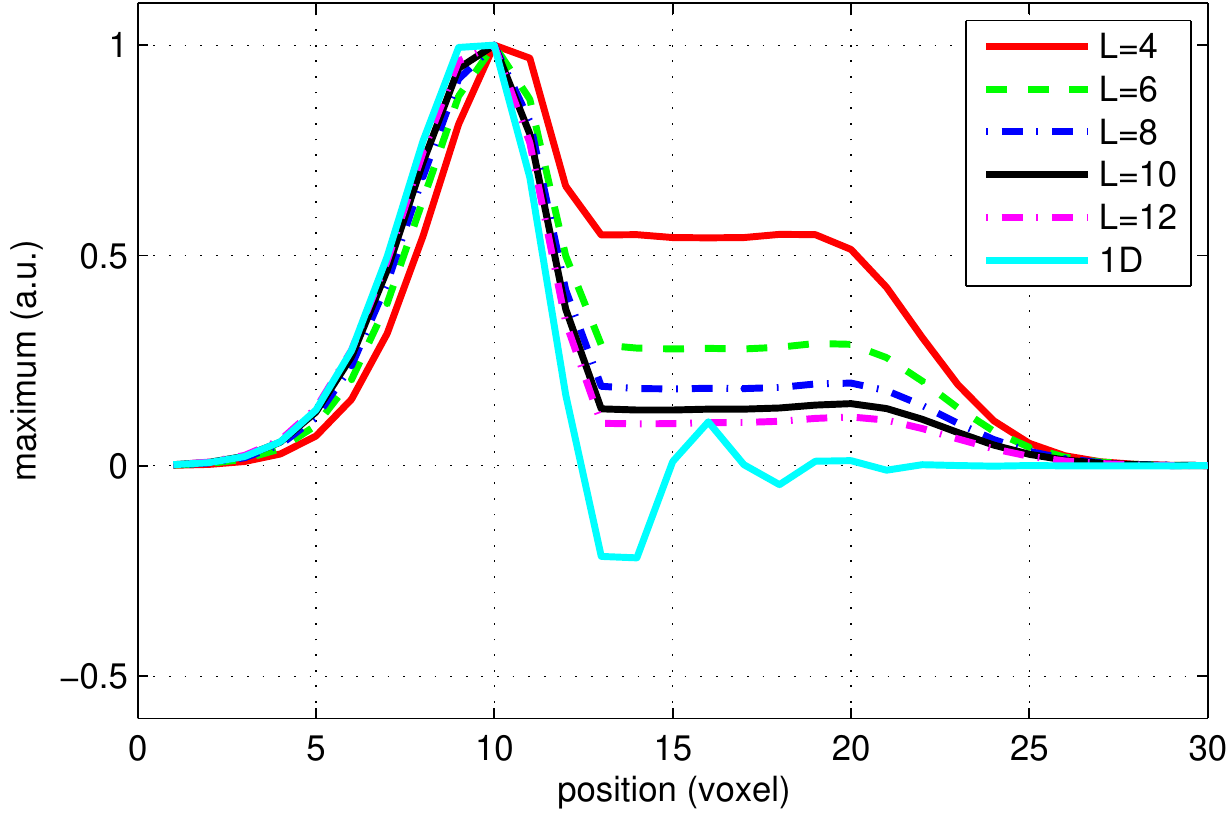}
\includegraphics[width=8cm]{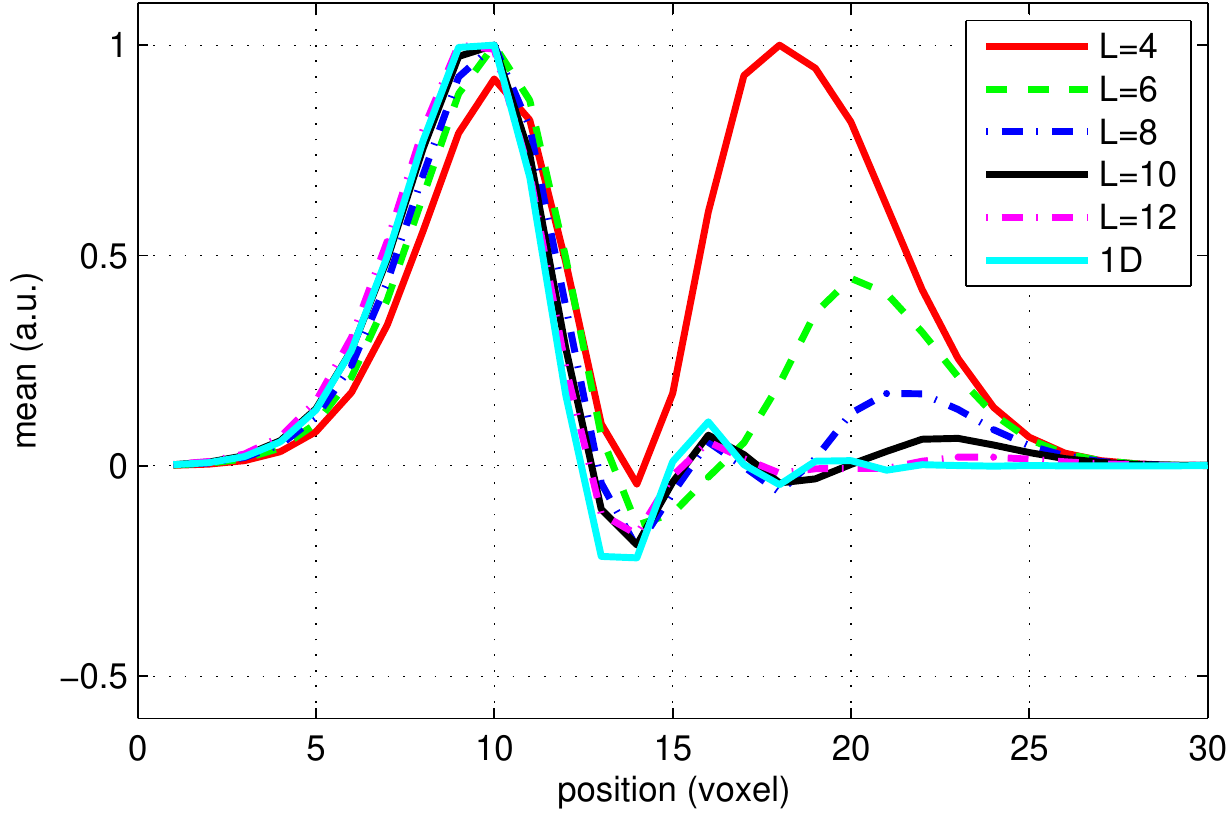}
\caption{The kernel $\Tk{0}$ was iterated $150$ times with a step-width of $0.05$ resulting in a translation of $7.5$ voxels. The initial condition was 
set to $\phi_0(\mv r,\mv n) = e^{-|\mv r|^2/2} \delta_{\mv n_z}(\mv n)$. 
On the left the maximum  $\max_{\mv n} \phi(z \mv n_z, \mv n)$ along the z-axis is plotted. On the right the mean-component $j=0$, i.e. $\mv f^0(z\mv n_z)$, 
is visualized. } \label{fig:translation_test2}
\end{figure}

\begin{figure}
\centering
\includegraphics[width=14cm]{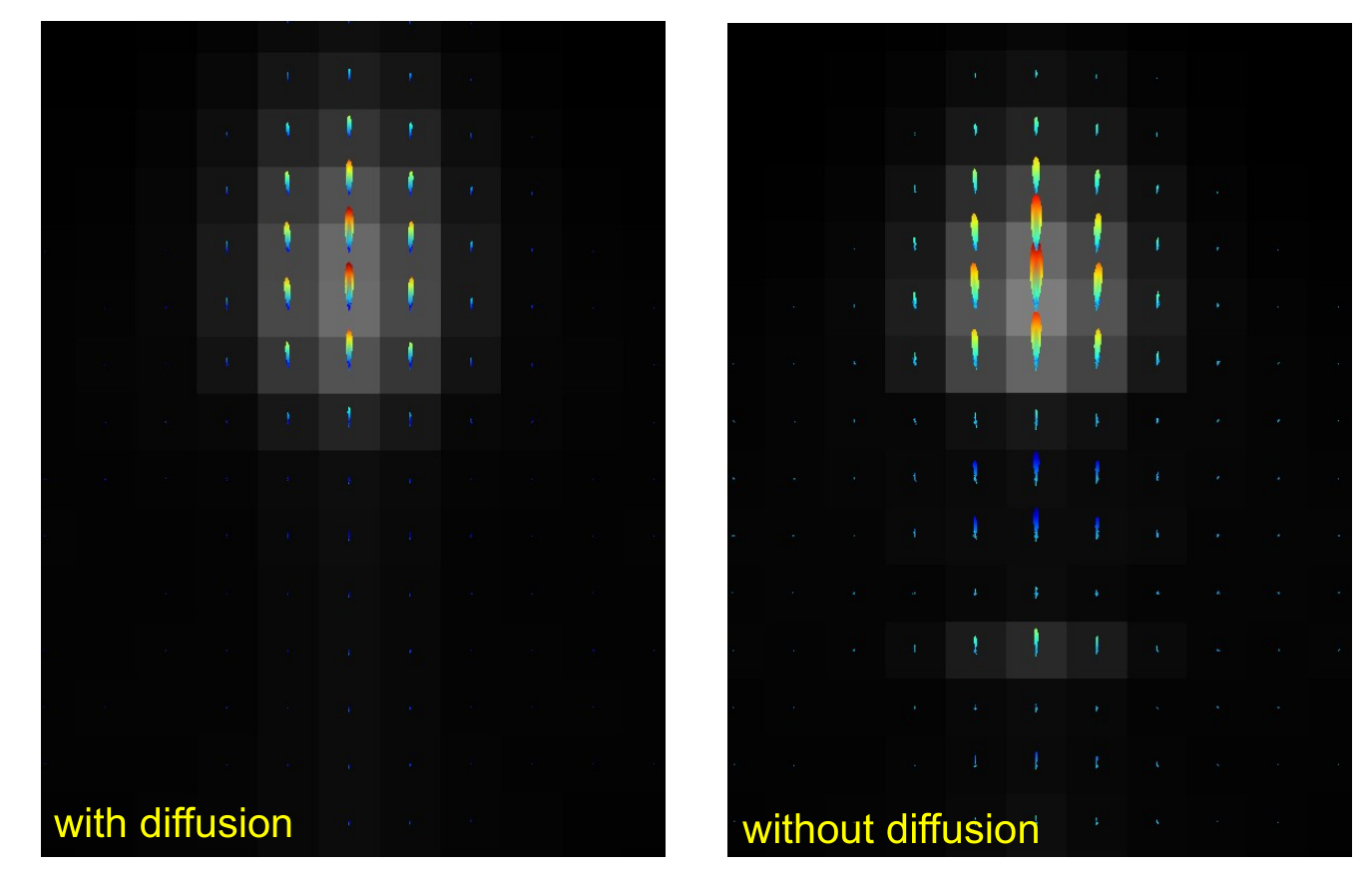}
\caption{The same setting like in \ref{fig:translation_test1} and \ref{fig:translation_test2} but the orientation distributions are visualized by glyphs.
The spherical harmonic cutoff is set to $L=12$. The underlying gray-values indicate the maximum of the orientation distribution.} \label{fig:translation_test3}
\end{figure}

In a first small experiment, we investigate the behavior of the simple forward integration of $\Tk{0}$ in detail. Additionally, we look at the forward integration of  $\Tk{0} + 0.1\Tk{0}^2$, i.e. a small 
diffusion term is added to make the results more stable. In particular, 
we used $\phi_0(\mv r,\mv n) = e^{-|\mv r|^2/2} \delta_{\mv n_z}(\mv n)$ as initial condition and iterated the integration $150$ times with a step-width of $0.05$ resulting in a translation 
of $t=150\cdot 0.05 = 7.5$ voxels. To analyze the result we have a look at the maximum of the orientation distribution along the z-axis, namely 
$\max_{\mv n} \phi(z \mv n_z, \mv n)$ and the $j=0$ component $\mv f^0(z\mv n_z)$, which is the mean along the orientation coordinate. For reference, we made a simple 1D experiment without any
orientation involved. It was just integrated $\phi(z) \mapsto \phi(z) + 0.05\cdot \partial_z \phi(z)$ with the same first-order central difference approximation of $\partial_z$. 
In Figure \ref{fig:translation_test1} the results for different cutoffs $L\geq j$ and the 1D experiment are shown. The 1D reference experiment shows the typical oscillations. Similar oscillations 
can also be observed for $j=0$ component. In particular, for low cutoffs the oscillations are pretty heavy and converge for large $L$ towards the 1D reference. For the maximum of the orientation 
distribution the behavior is different. Instead of oscillations a tail is left over, while certain peaks appear where the oscillations of $j=0$ component dominate. 
Looking at Figure \ref{fig:translation_test2} where the small additional diffusion term is considered, one can observe less oscillations and a more regular behavior. 
In Figure \ref{fig:translation_test3} we show the resulting orientation fields for $L=12$ in glyph representation with and without diffusion. Once again 
one can see the spike for the version without diffusion and the more regular behavior for the diffusion-regularized version.

%
%
%
%
%
%
%

\newcommand{\mo}{h}
\newcommand{\Mo}{\mv H}
\subsection{Spatially Regularized Spherical Deconvolution}
Magnetic resonance imaging (MRI) has the  potential to visualize non-invasively the fibrous structure of the human brain
 white matter \cite{Jones2010}. 
Based on the  directional dependency of water diffusivity in fibrous tissue it is possible to reveal underlying connectivity information. 
The accurate and reliable processing and estimation of fiber orientation distributions is 
a major prerequisite for the processing of such data. 
There are numerous methods for estimating orientation distributions on the basis of the diffusion-weighted MR-signal.
We will focus spherical deconvolution \cite{tournier2004}, which is 
one way to estimate the so called fiber orientation distributions (FOD) on the basis of the diffusion-weighted MR-signal. The idea is based on 
a model-driven deconvolution scheme to turn the diffusion weighted MR-signal into a FOD. 
The goal is to find a FOD $f$ such that
\begin{equation} 
J(f) = \iint\limits_{\Real^3 \times S_2} \left| (\Mo f)(\mv x, \mv n) - S(\mv x,\mv n) \right|^2 \ d\mv x d\mv n \label{eq:deconv}
\end{equation}
is minimized. Here $S$ denotes a quantity derived from the MR-measurement, e.g. in standard q-ball imaging just the ratio $M(\mv x,\mv n) / M_0(\mv x)$,
where $M$ is the diffusion weighted image and $M_0$ the measurement without diffusion weighting.
The operator $\mv H$ denotes the spherical convolution with the so called fiber response function:
\begin{equation}
(\Mo f)(\mv x, \mv n) = \int_{S_2} \mo(\mv n \cdot \mv n') f(\mv x,\mv n') d\mv n',
\end{equation}
The fiber response function is typically chosen to be a function of the form $\mo(t) = \exp(-\lambda t^2)$, or  is estimated from the measurement itself. 
The main problem is that $\mv H$ is practically not invertible. One way to solve the problem is to introduce a non-negativity constraint \cite{Tournier2007}, 
but in case of low-quality data the problem is still
hard to invert. Therefore, regularization techniques \cite{reisert2011,duits:leftinv,Goh2009} were introduced, which make the solution unique and stable. 
An additional term $R(f)$ is added to the cost function as:
\begin{equation}
J_\text{reg}(f) = J(f) + R(f) \label{eq:objective} \nonumber
\end{equation}
The fiber continuity/contour enhancement kernel provides such an additional cost function that perfectly fits to the nature 
of fiber orientation distributions:
\begin{equation}
R_\text{FC}(f) =\lambda \iint\limits_{\Real^3 \times S_2} (\mv n\cdot\nabla f)^2\ d\mv x d\mv n. = \lambda  \iint\limits_{\Real^3 \times S_2} (\Tk{0} f)^2\ d\mv x d\mv n.
\end{equation}
We can expect from this kernel that it prevents `arbitrary' smoothing and that it will preserve and emphasize the fibrous nature of the data. 
As the regularizer includes the spatial neighborhood of each voxel, we have
to be careful at the transition area between gray and white matter: a boundary condition is needed.
Instead of a hard boundary condition we decided to keep the FOD small in the background area (non white matter). Thus, 
we have an additional term in the cost function that suppresses the FOD in the background,
\begin{equation}
R_\text{mask}(f) = \lambda_\text{mask} \iint\limits_{\Real^3 \times S_2} f(\mv x,\mv n)^2\ (1-w_\text{m}(\mv x))\ d\mv x d\mv n,
\end{equation}
where $w_\text{m}(\mv x)$ is a white matter mask. 
That is  $\lambda_\text{mask}$ acts as a suppression strength $\lambda_\text{mask}$ for gray-matter. We found a value of $1$ to be a good value.

The choice of regularization strength is a crucial issue. We found that too strong regularization emphasizes discretization artifacts, which 
is shown by a violation of the rotation invariance, that is, the results depend on the absolute directions of the bundles, with an unpredictable dependence on the sphere 
discretization and the underlying voxel grid. On the other hand too low values lead to less stable results. We found by a simple visual inspection of a simulated crossing 
a value of $\lambda= 0.005$ to be a good trade-off. 

\subsubsection{Optimization and Implementation}
In order to find the optimum of $J_\text{reg}$, we have to compute the variation of the objective and set it to zero for
the necessary condition for the optimum:
\begin{equation}
\delta_f J_\text{reg} = \Mo^\T \Mo f  - \Mo^\T S   - \lambda \Tk{0}^2 f - \lambda_\text{mask} (1-w_m)f = 0
\end{equation}
which leads to 
\begin{equation}
(\Mo^\T \Mo  - \lambda \Tk{0}^2 - \lambda_\text{mask} (1-w_m))f =  \Mo^\T S
\end{equation}
To solve this linear equation we employed an ordinary conjugate gradients scheme.
The implementation of the operator on the left-hand side in terms of spherical harmonics is based on 
equation \eqref{eq:T2res}. The operator $\Mo$ is also quite simple to implement in terms of spherical harmonics, because it is a diagonal matrix.
The elements on the diagonal are just $\int_{-1}^1 h(t) P_j(t) dt$. If we want to implement this with a discretized sphere the operator $\Mo$ 
has to approximated by interpolation. We just used spherical harmonics to do this interpolation, in the same 
way like in \cite{duits:leftinv,duits:diffconv20011} the spherical Laplace-Beltrami operator was approximated. The operator $\Tk{0}^2$ is discretized very similar
to our spherical harmonics implementation with finite differences (comparable to \cite{duits:diffconv20011}, Appendix F). Throughout this experiment 
the sphere was discretized with $512$ direction, which 
were determined such that the energy of the configuration of 512 electrons on the surface of a sphere is minimized, where the electrons repel 
each other with a force given by Coulomb's law.

For both methods the conjugate gradients algorithm was iterated $100$ times, which was enough for convergence.

\begin{figure}
\centering
\includegraphics[width=16cm,height=11cm]{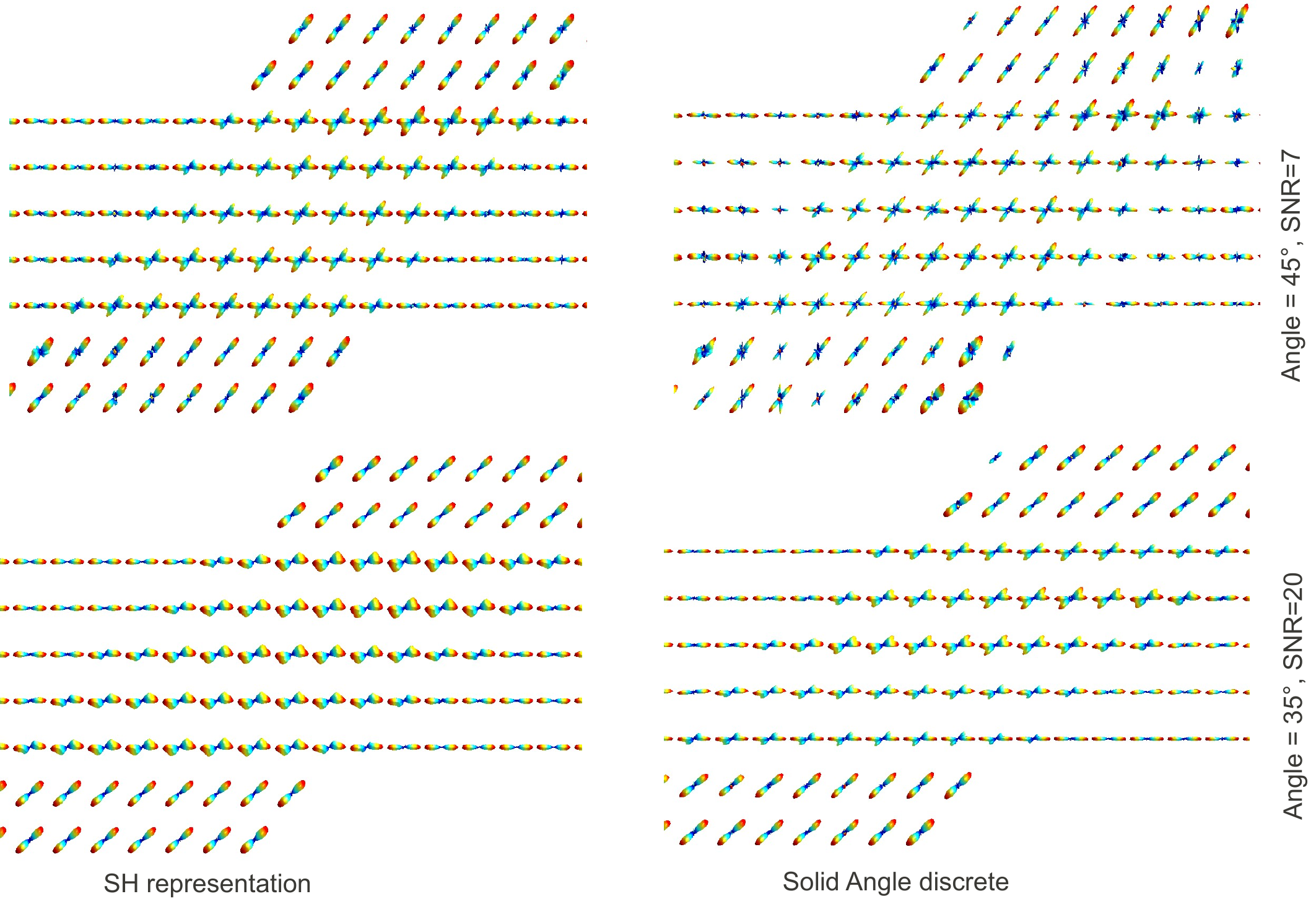}
\caption{Two simulated crossing situations to compare SH representation with the angular discrete version. The simulation was performed with $64$ gradient directions at a 
b-value of  $1000 s/mm^2$ and with a diffusion coefficient of $10^{-3} mm^2/s$} \label{fig:crossing}
\end{figure}

\subsubsection{Experimental Setup}
The goal of the experiment is to compare the proposed spherical harmonic representations with the discrete angular representation of the fiber orientation distributions.
Therefore, we simulated the MR-signal with 64 gradient directions of a crossing region. 
The signal was simulated by using the standard exponential model $S_{\mv n_\text{fib}}(\mv n) = e^{-bD(\mv n \cdot \mv n_\text{fib} )^2}$, that is, no diffusion perpendicular 
to a fiber is assumed. We have chosen  $bD = 1$, which emulates a b-value of $1000 s/mm^2$ and a typical diffusion coefficient for the human brain. 
The generated signal was distorted by Rician noise $S_\text{noisy} = \sqrt{(S + n_\text{real})^2 + n_\text{imag}^2}$, where 
$n_\text{real}$ and $ n_\text{imag}$ are normally distributed real numbers with standard deviation $\sigma$. 
The signal-to-noise ratio is defined as $\text{SNR} = 1/\sigma$, that is, the SNR is calculated with respect to the $b=0$ measurement.
The crossing was created on a  $24\times 24$ voxel grid, where the tracts of the crossing are on average 5 voxels thick.
To get an impression look at Figure \ref{fig:crossing}.

To measure the performance of the deconvolution method, the local maxima of the estimated FODs are extracted
and compared to the ground truth direction. To efficiently search for a local maximum we followed the this procedure:
We made a Voronoi tessellation of the sphere\footnote{in the case of the SH-representation the orientation distribution was sampled with the same 
directions like the angular discrete algorithm uses}
and search for those direction whose value is above all its neighbor with respect to the tessellation.
To accurately determine the underlying continuous direction, a quadratic form is fitted to the neighborhood of the putative local maxima and the maximum 
of this quadratic form is used as detected direction. A ground truth direction is said to detected, when it is in a range of $10^\circ$ from a detected direction.
To measure the performance we used precision, recall and the f-score\footnote{Let TP be the number of successfully found ground truth directions,
let FP be the number of detections that are not in a range of $10$ degree to a ground truth direction, and let FN the number of ground truth direction that 
are not detected, then $\text{precision} =  \text{TP} / ( \text{TP} + \text{FP})$ and $\text{recall} =  \text{TP} / ( \text{TP} + \text{FN})$ and 
$\text{f-score} =  2 \ \text{precision}\cdot \text{recall} / ( \text{precision} + \text{recall} )$.}.
To generate the performance measures the simulations were repeated 100 times.

\begin{figure}
\centering
\rotatebox{90}{\hspace{1.5cm} $\alpha = 0^\circ$}
\includegraphics[width=5cm]{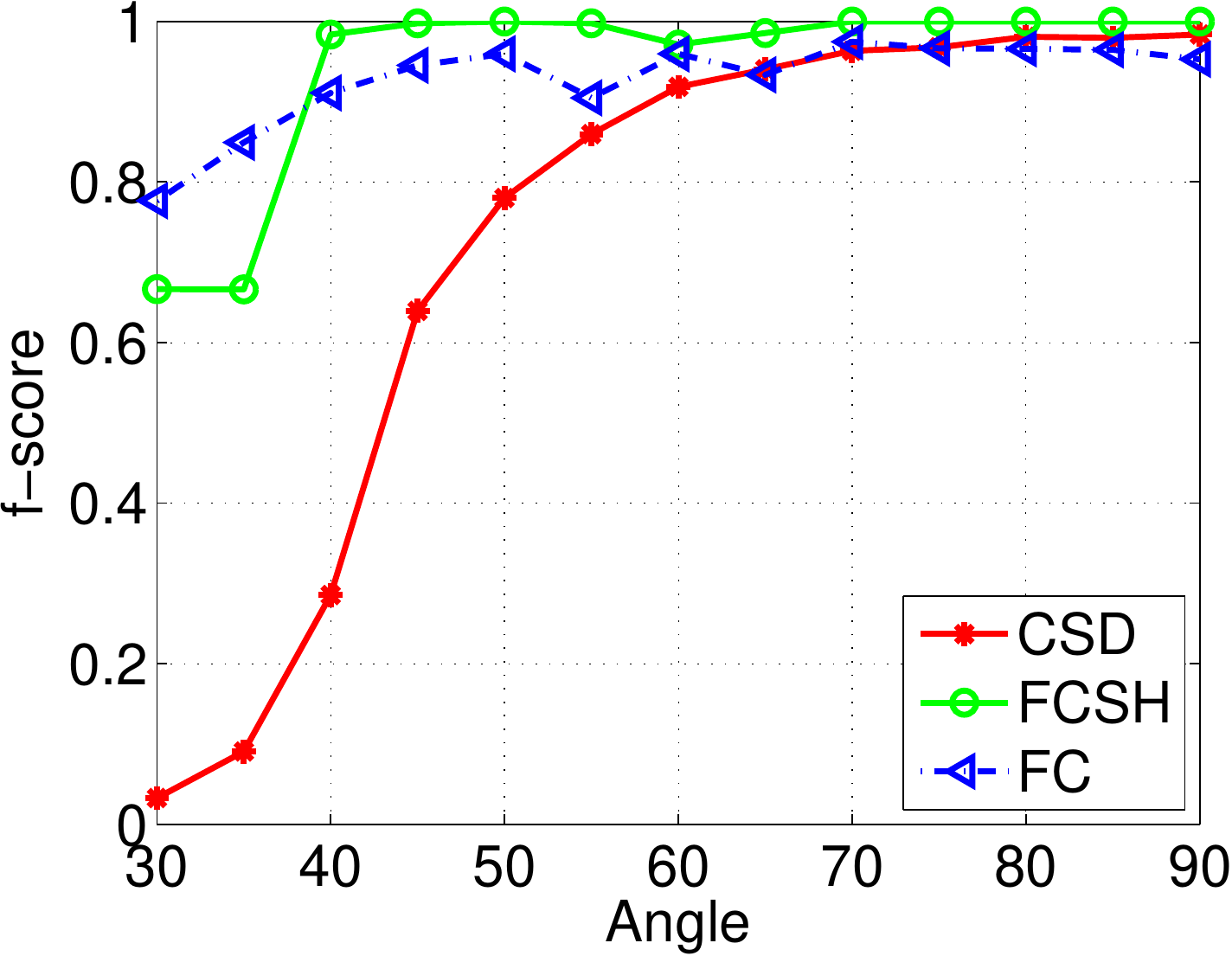}
\includegraphics[width=5cm]{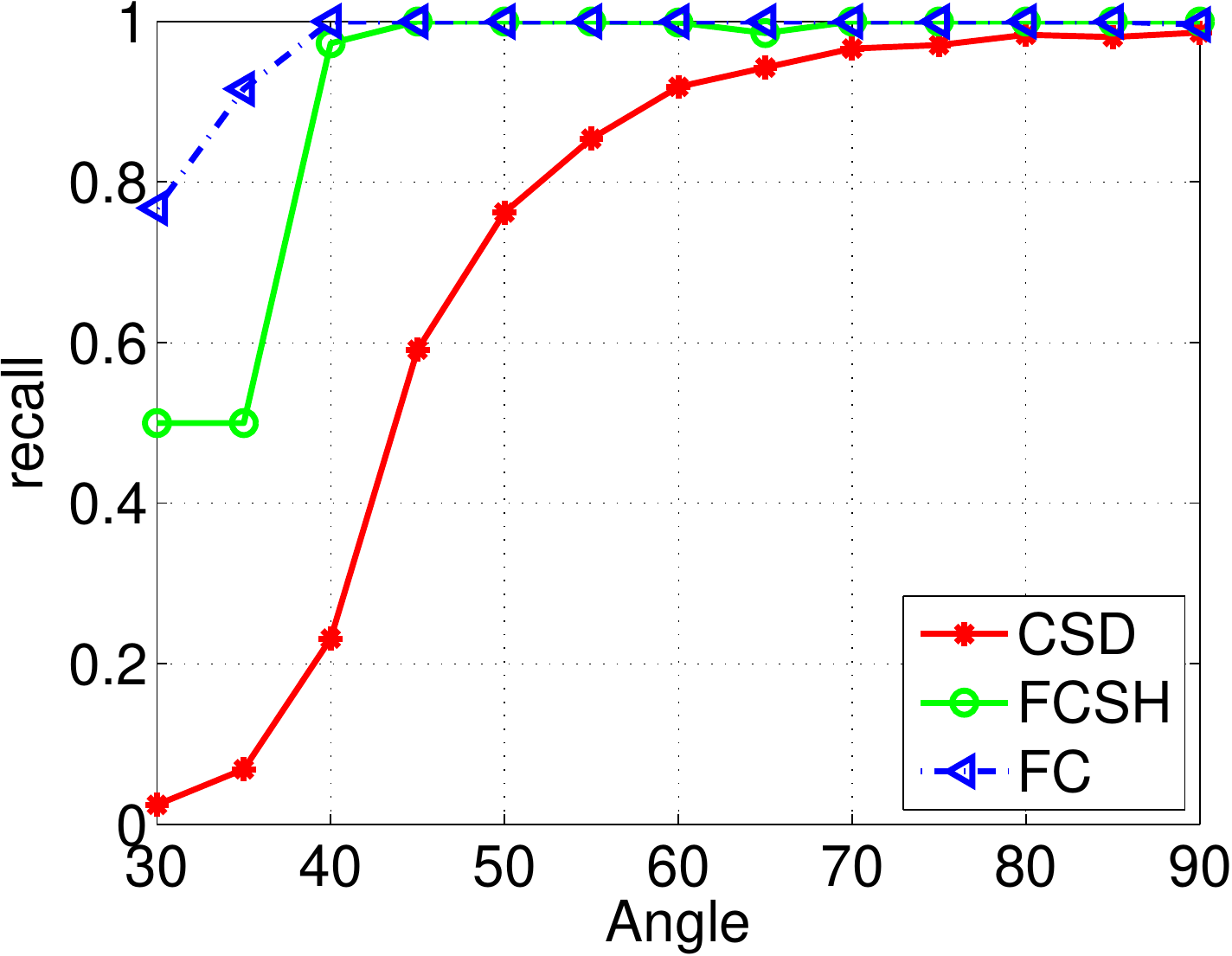}
\includegraphics[width=5cm]{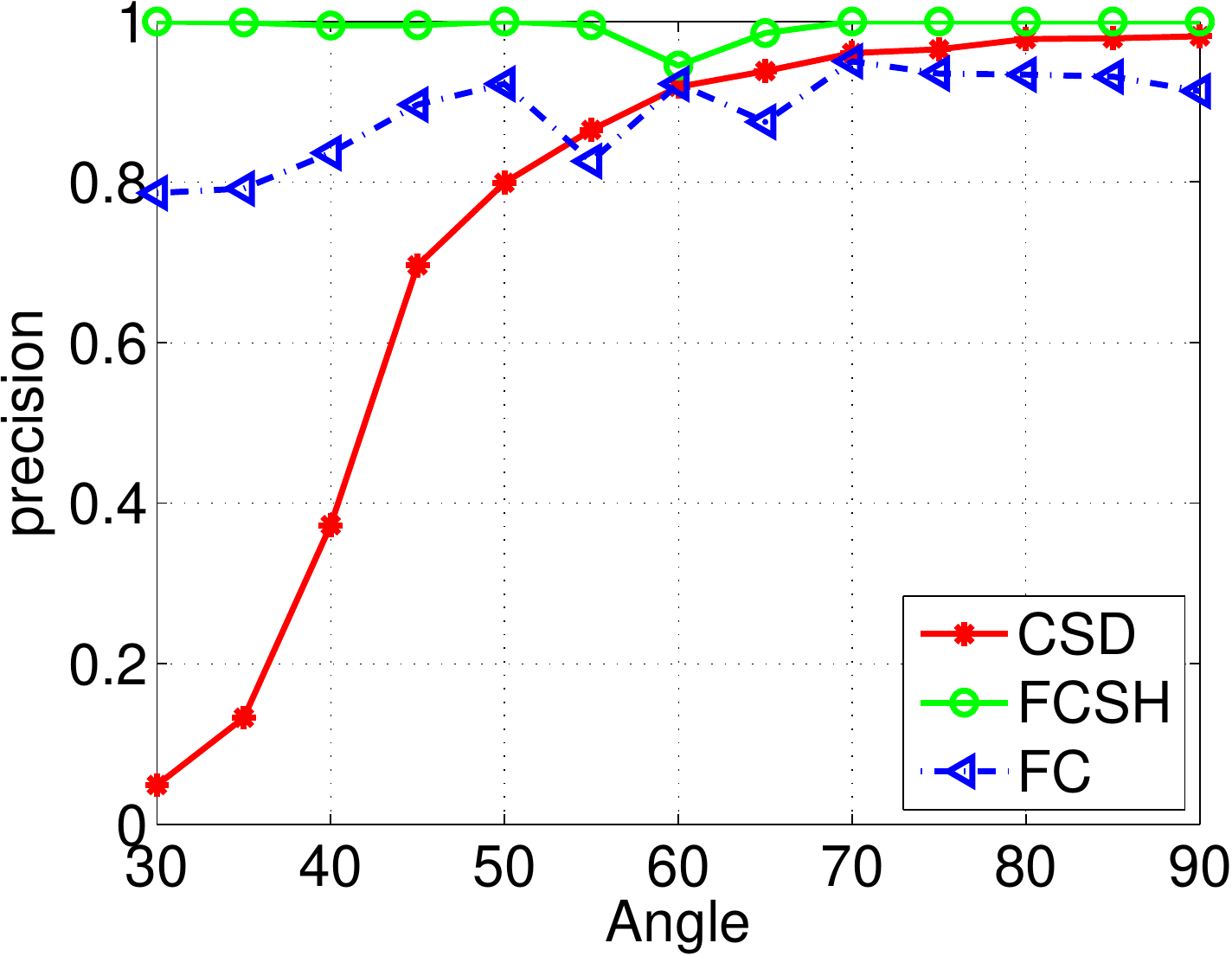} \\
\rotatebox{90}{\hspace{1.5cm} $\alpha = 15^\circ$}
\includegraphics[width=5cm]{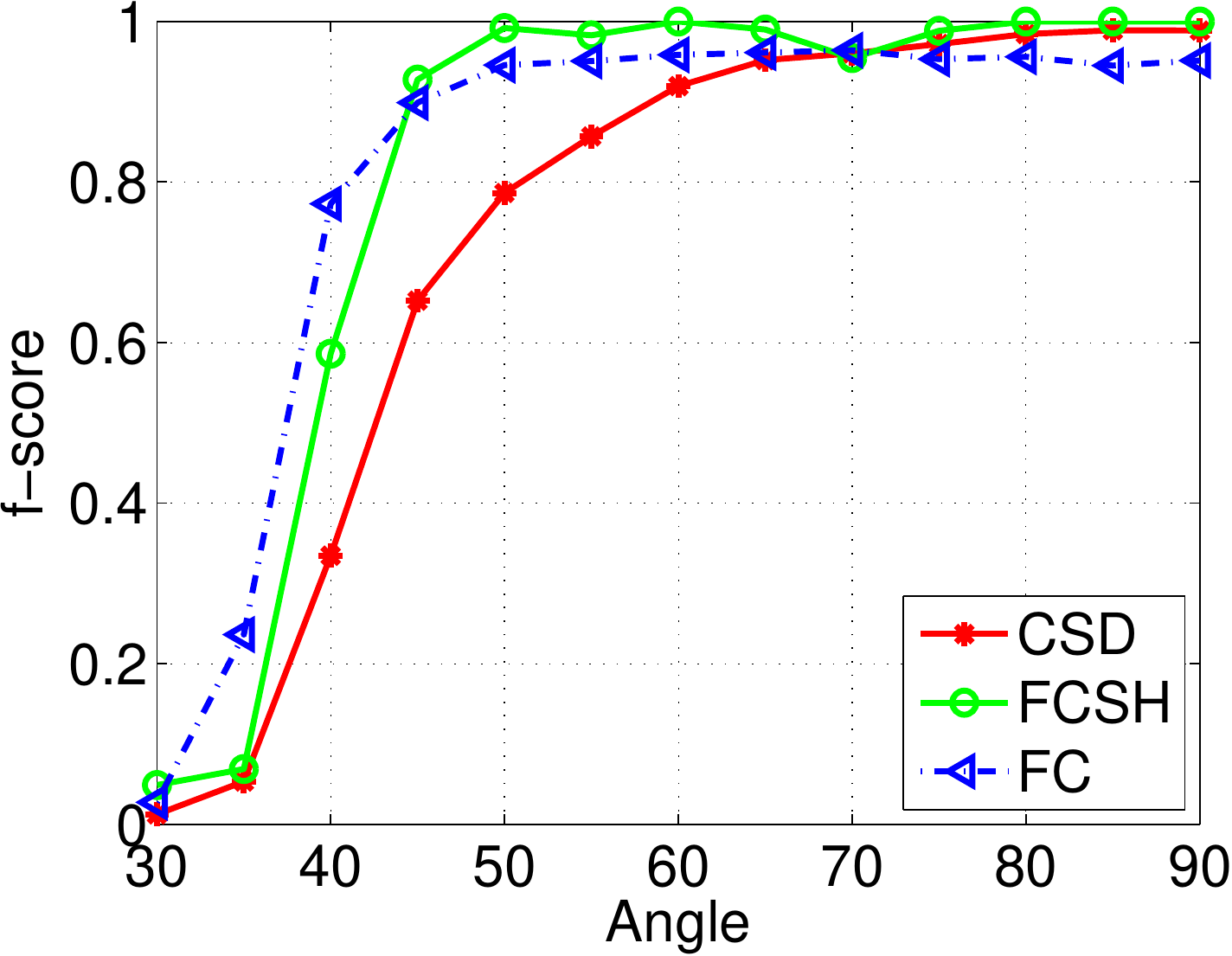}
\includegraphics[width=5cm]{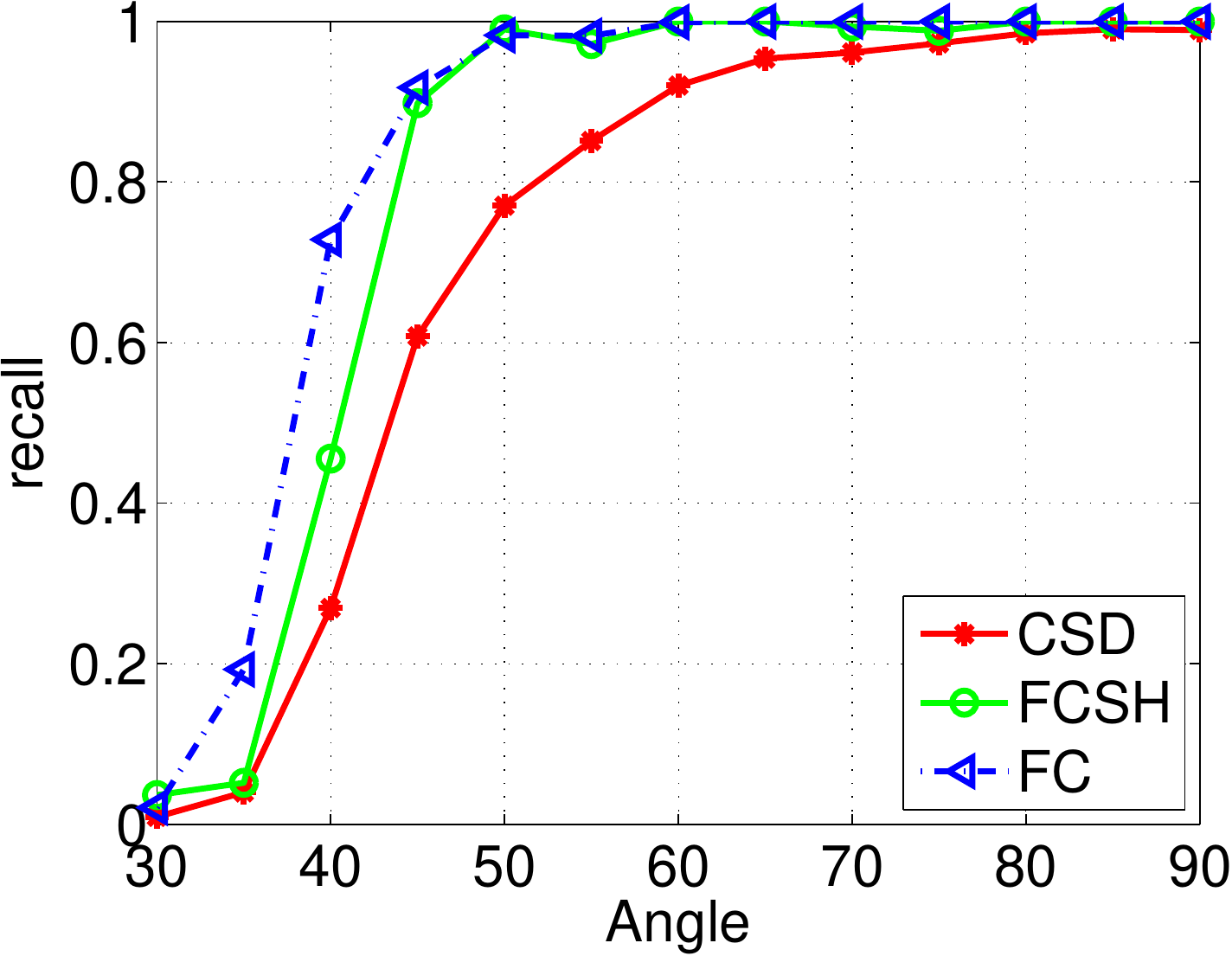}
\includegraphics[width=5cm]{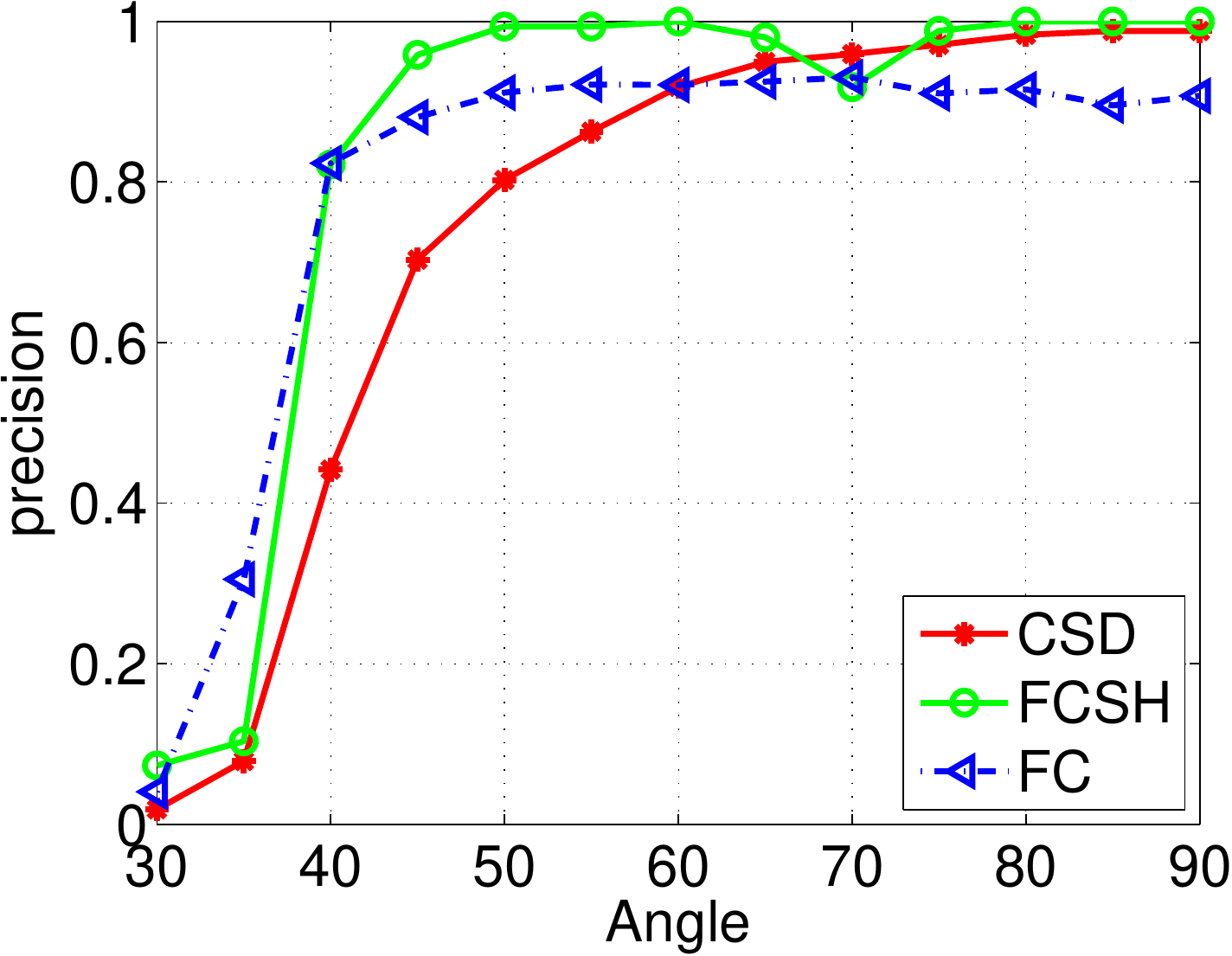} \\
\rotatebox{90}{\hspace{1.5cm} $\alpha = 30^\circ$}
\includegraphics[width=5cm]{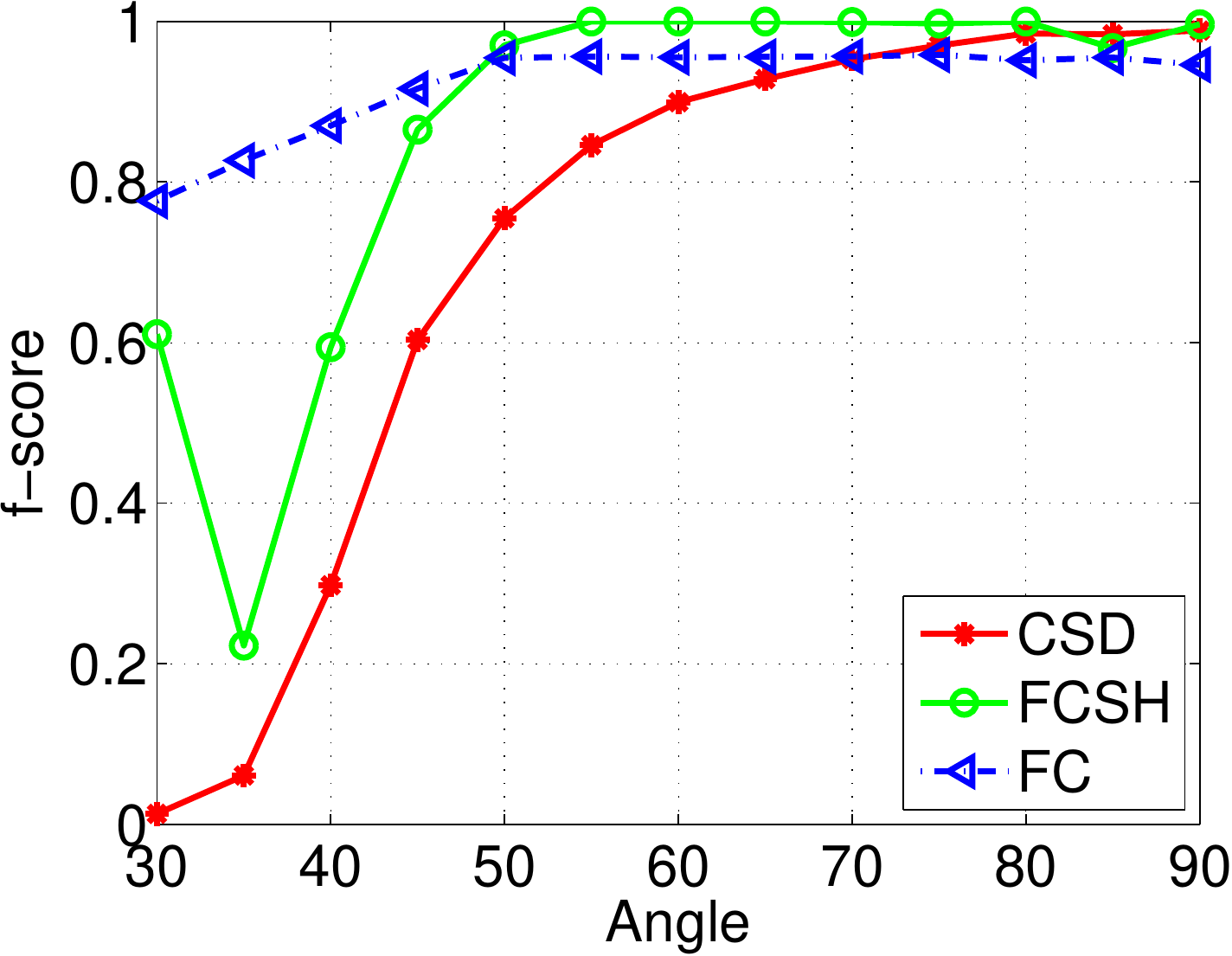}
\includegraphics[width=5cm]{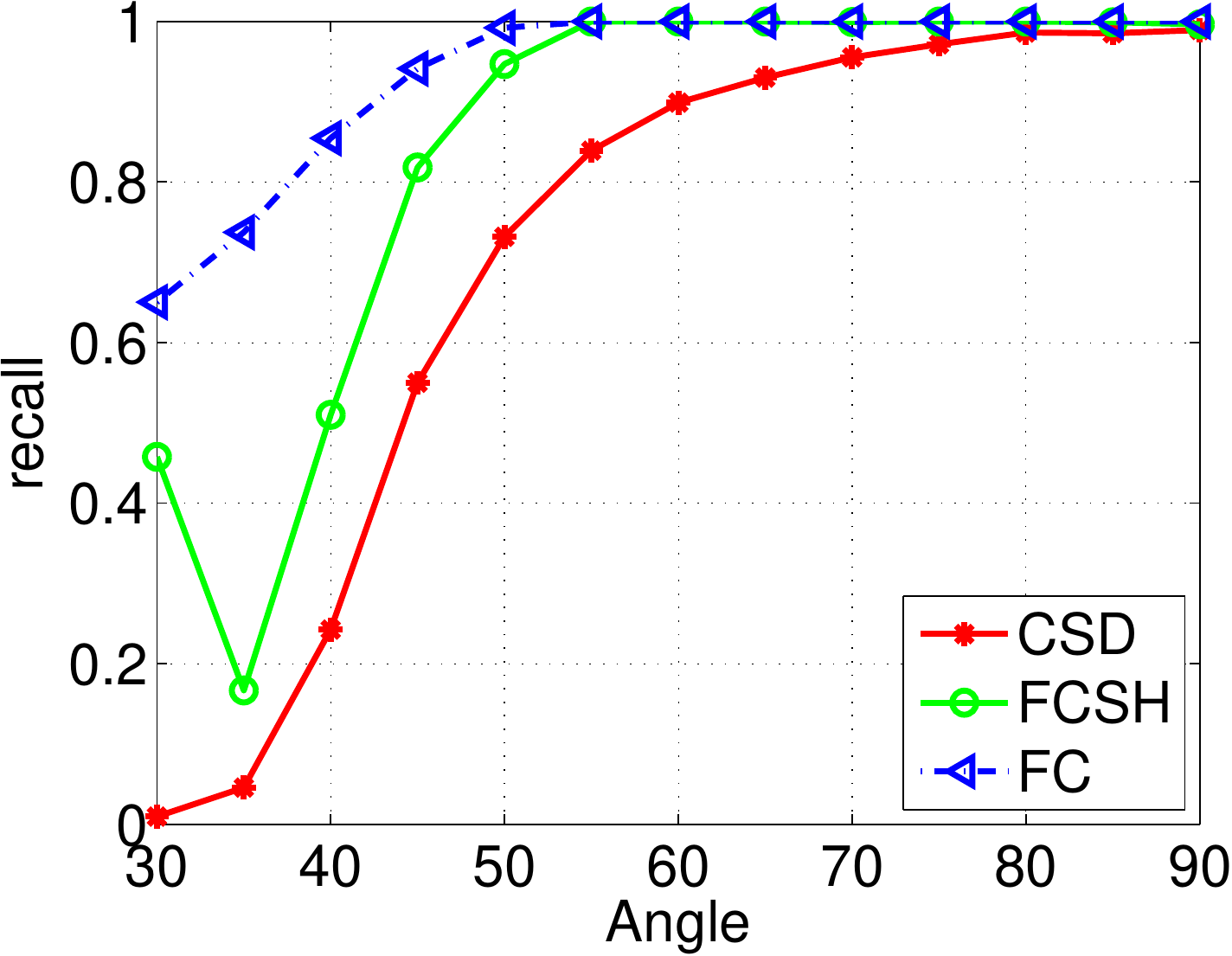}
\includegraphics[width=5cm]{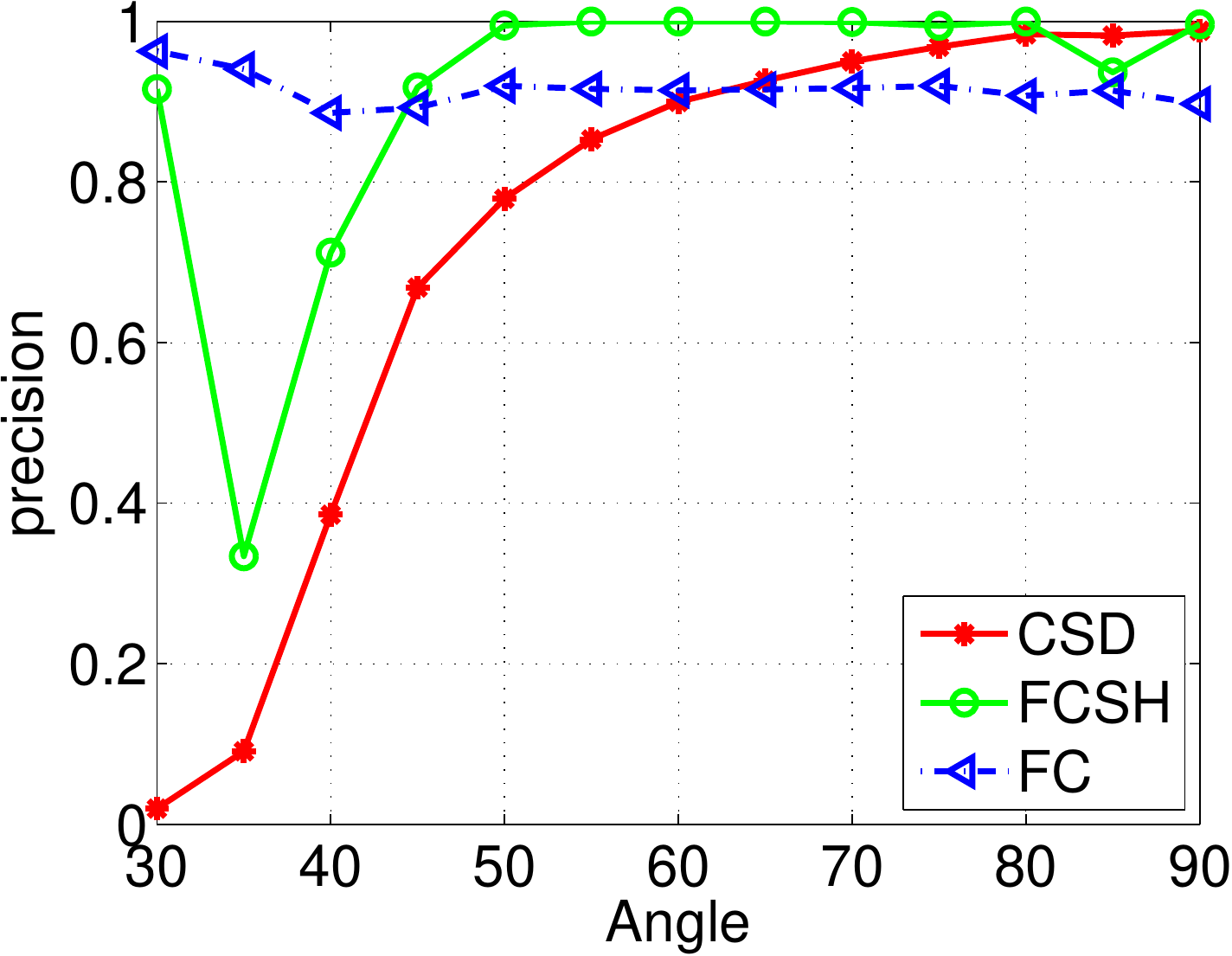} \\
\rotatebox{90}{\hspace{1.5cm} $\alpha = 45^\circ$}
\includegraphics[width=5cm]{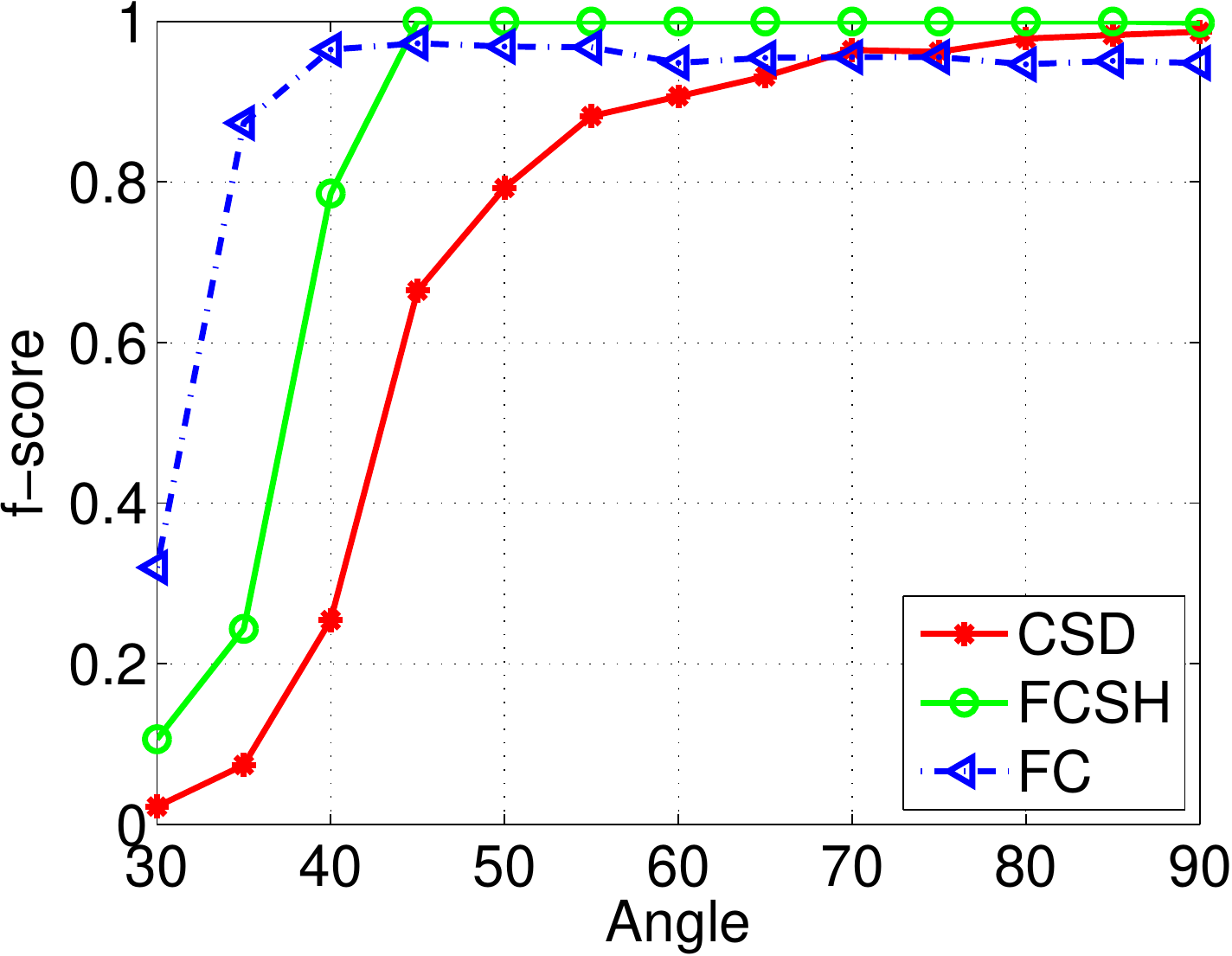}
\includegraphics[width=5cm]{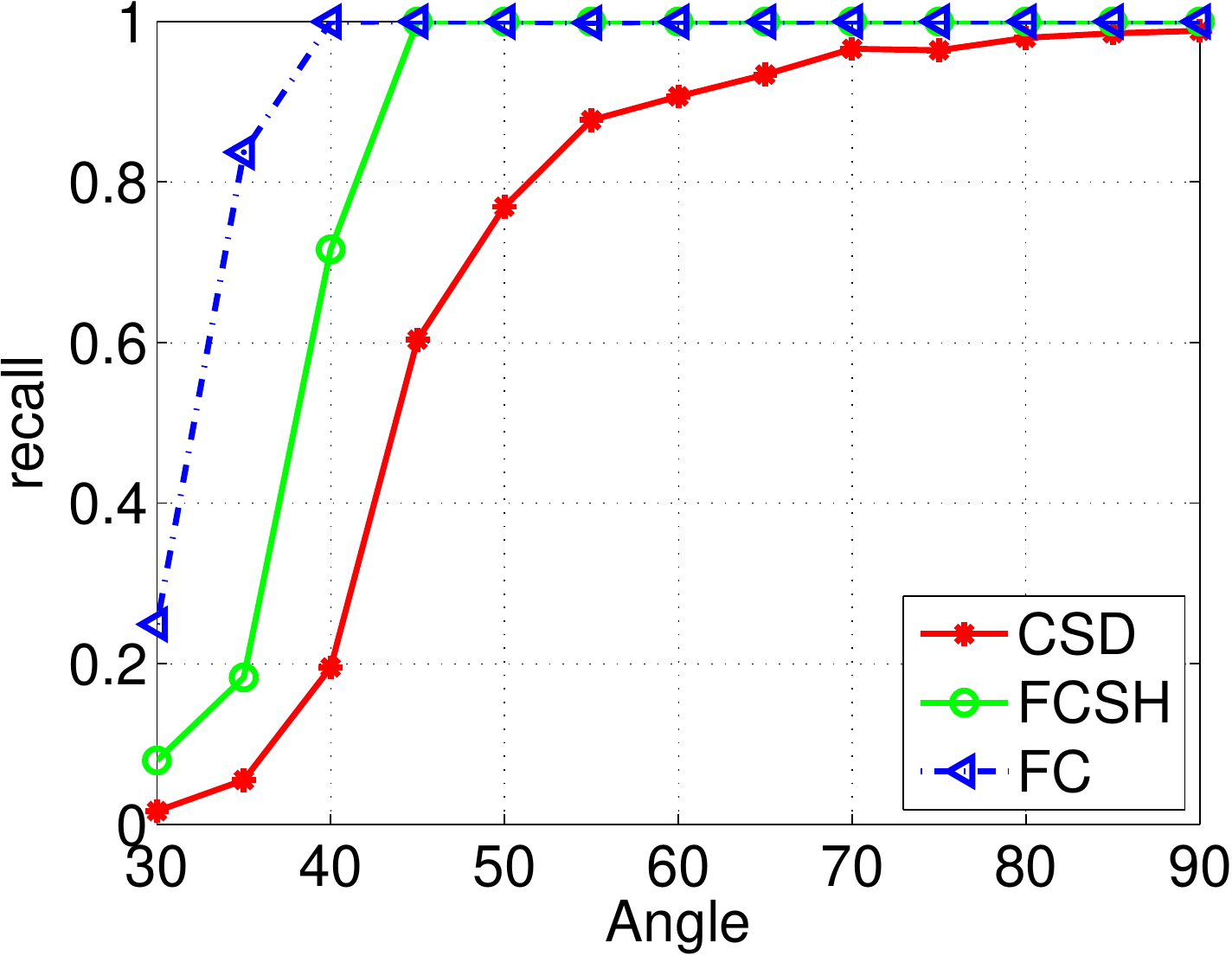}
\includegraphics[width=5cm]{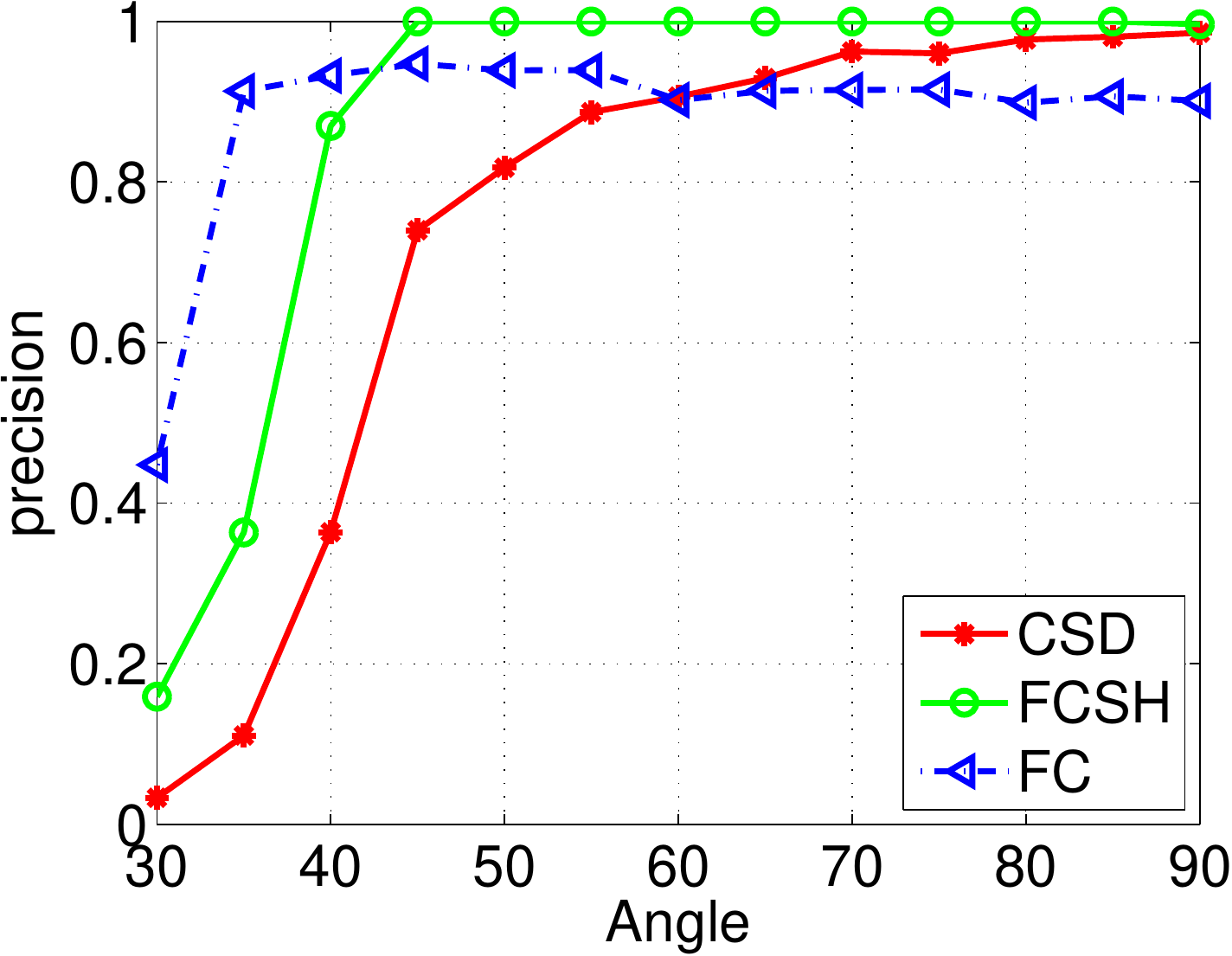} \\
\rotatebox{90}{\hspace{1.5cm} $\alpha = 60^\circ$}
\includegraphics[width=5cm]{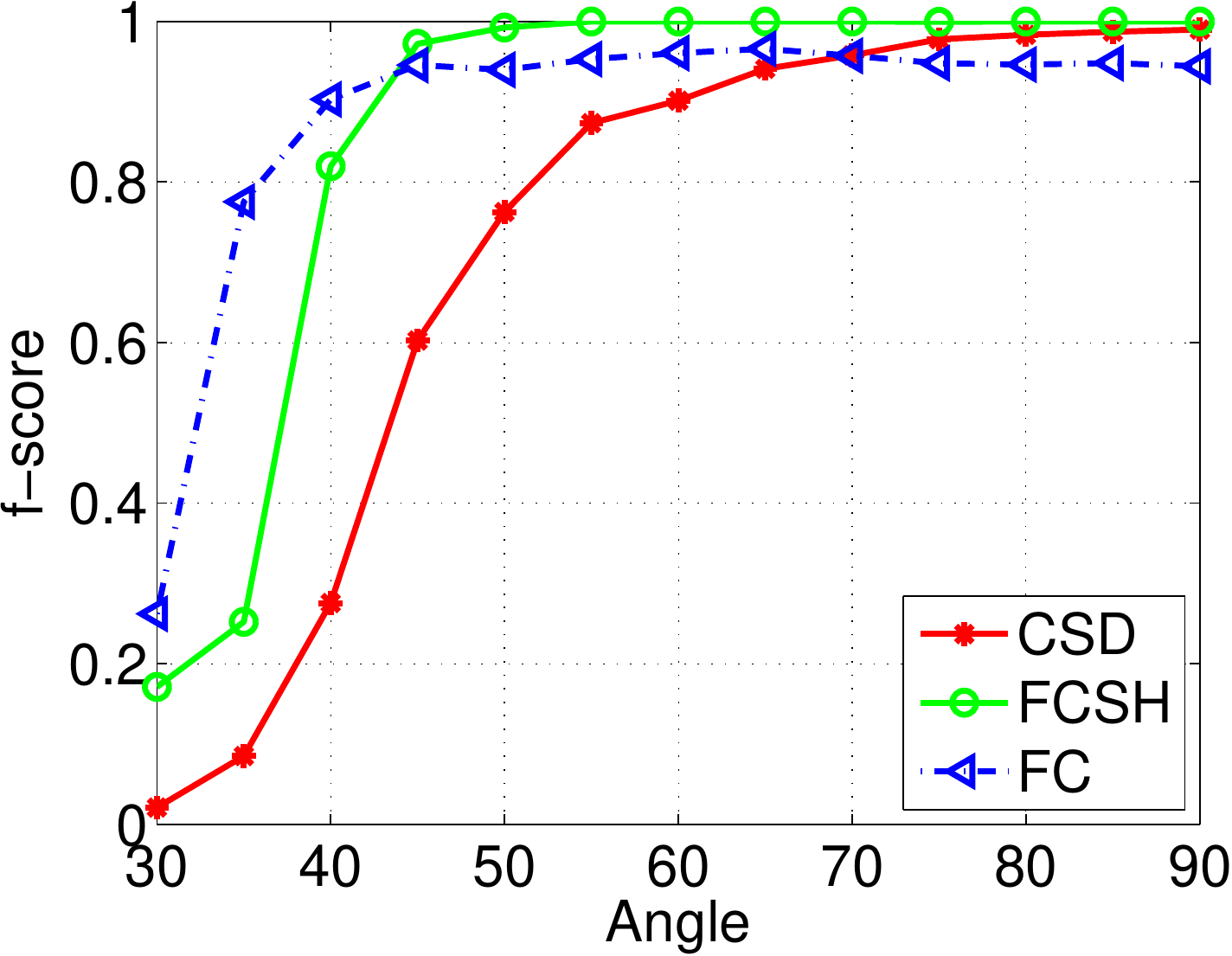}
\includegraphics[width=5cm]{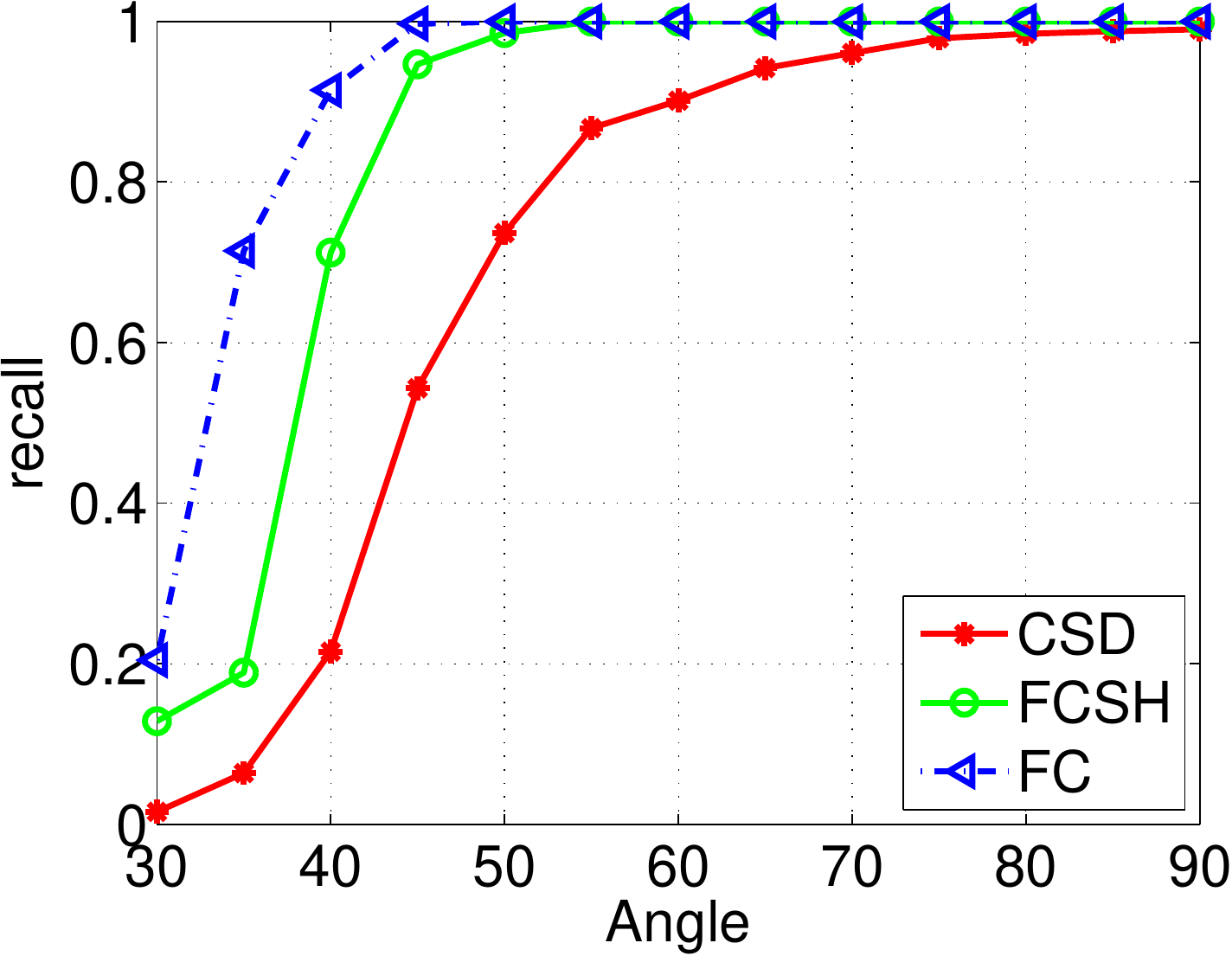}
\includegraphics[width=5cm]{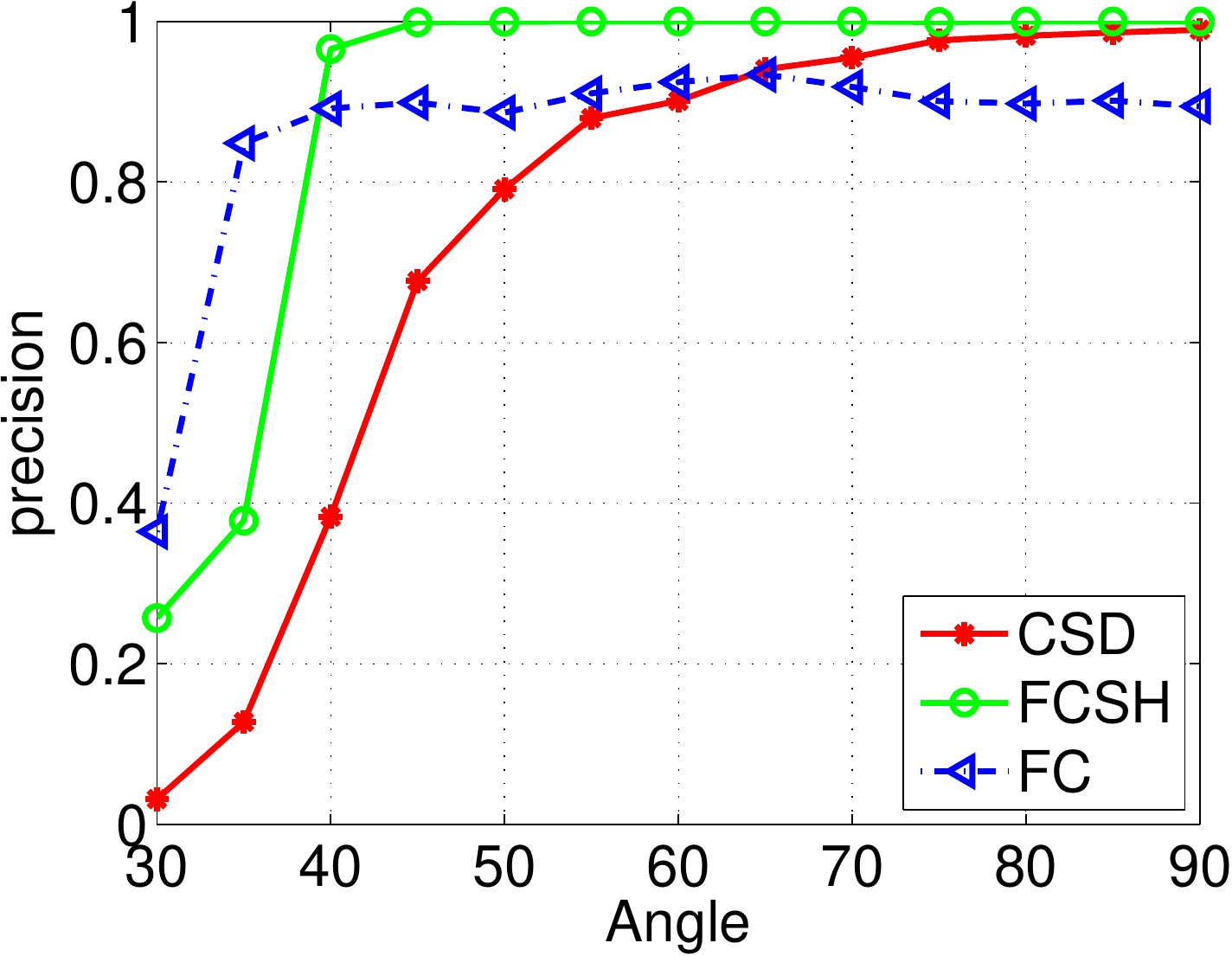} \\
\caption{Detection accuracies in terms of precision/recall and f-score for different crossing angles and absolute angles $\alpha$. 
The measurement was simulated at a $SNR=50$ with a value of $bD=1$ and $64$ gradient directions. 
} \label{fig:results}
\end{figure}

\begin{figure}
\centering
\includegraphics[width=17cm]{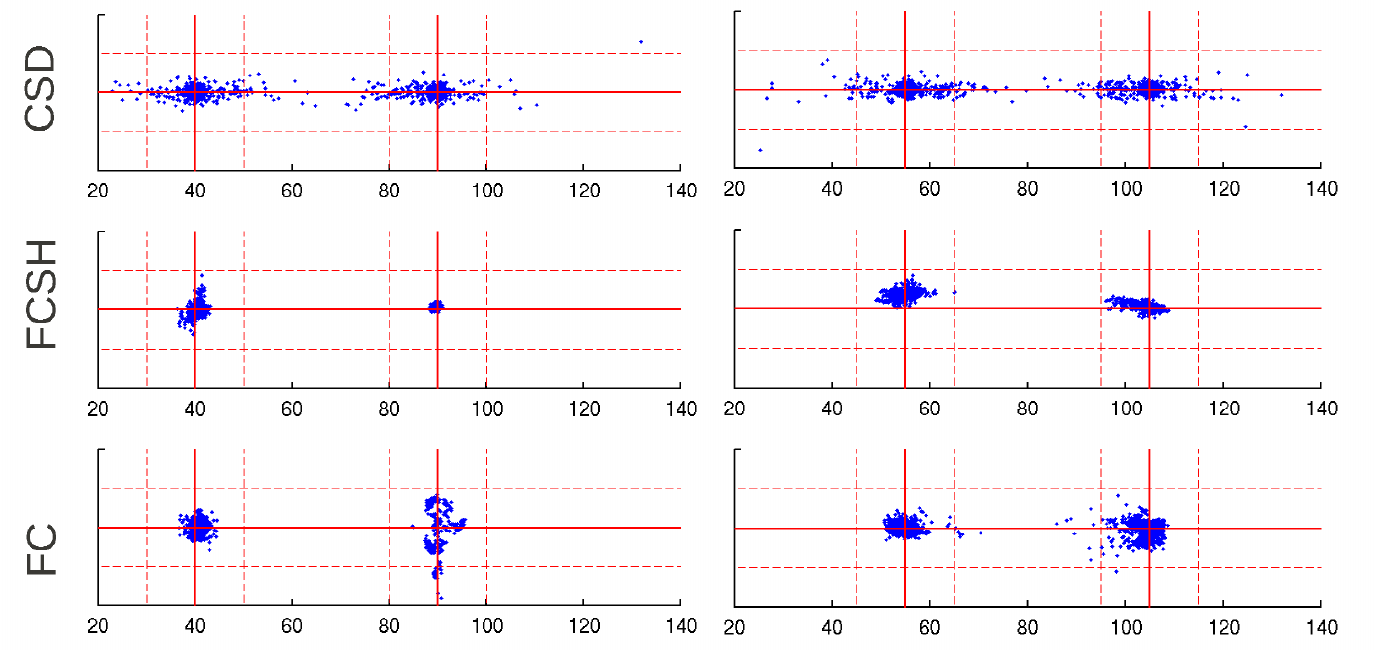}
\caption{Scatter plots in the $\phi$-$\theta$ plane for the crossing configuration with an crossing angle of $50^\circ$. 
The angle $\phi$ is plotted along the x-axis, $\theta$ along the y-axis of the scatter plots. The intersections of the thick red lines indicate 
the expected ground truth directions. The dotted lines indicate the $10^\circ$
 detection tolerance. On the left, the results for the configuration with an absolute angle of $\alpha=0$ are shown, on the right the results for 
$\alpha=15^\circ$ are given.} \label{fig:scatter}
\end{figure}

\subsubsection{Method Comparison and Discussion}

\begin{figure}
\centering
\includegraphics[width=16cm]{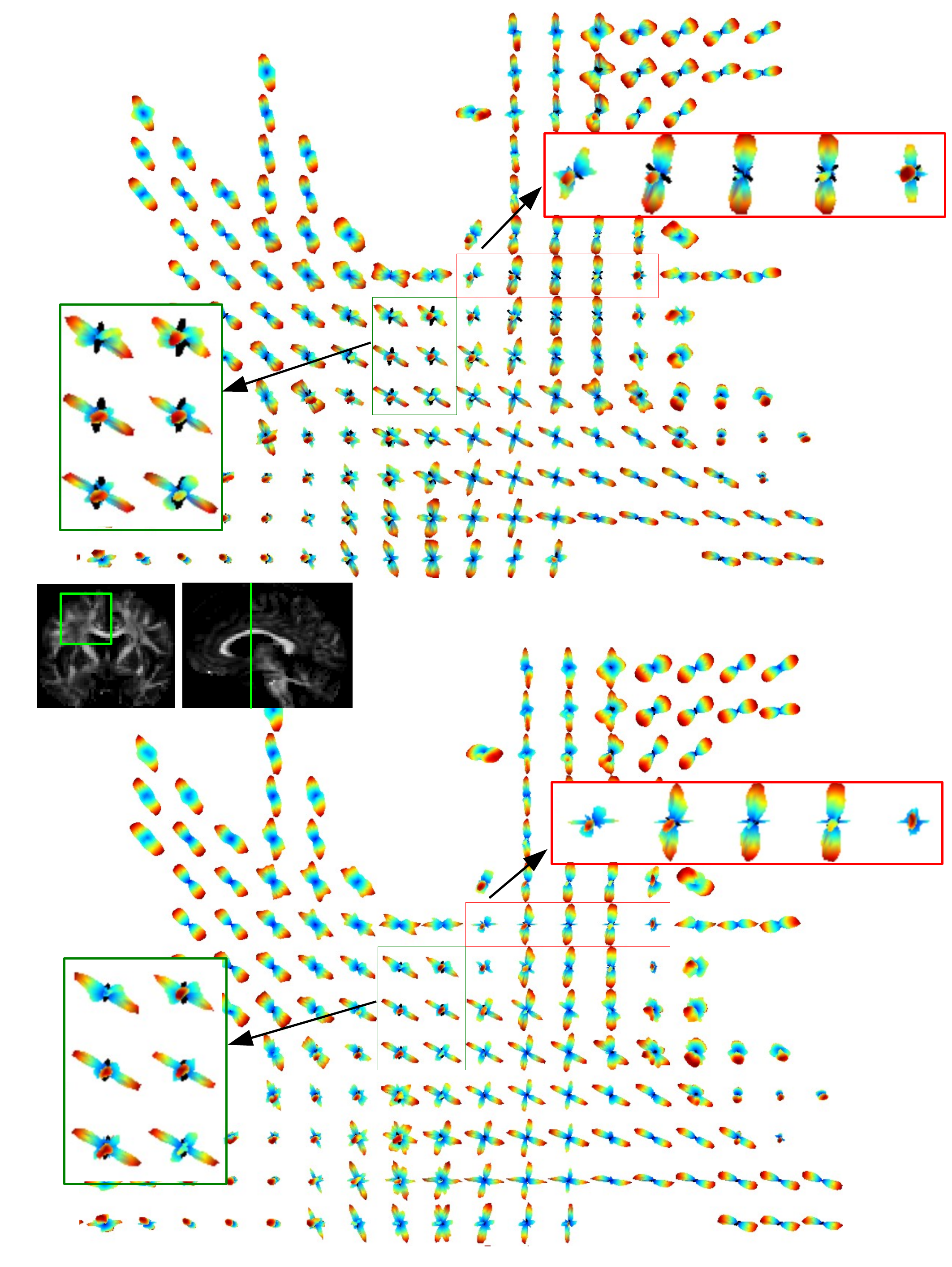}
\caption{A real world example of the human brain. A coronal section is shown. On top the results of the SH-based approach with $L=8$, on the bottom the discrete approach
with $512$ directions. } \label{fig:invivo_glyph}
\end{figure}

In Figure \ref{fig:crossing} we compare visually the spherical harmonic implementation for a cutoff of $L=8$ with the discrete version with $512$ directions on the sphere. 
We consider a crossing angle of $35^\circ$ and $45^\circ$. For both situations one fiber bundle direction was chosen along the underlying Cartesian coordinate axis. 
For the larger crossing, which is easier to resolve, we assumed a relatively low SNR of $7$. For this case one can 
see that the discrete version is more susceptible to noise than the spherical harmonic version, which is not astonishing due to the implicit regularization by the 
finite SH-cutoff of $L=8$. For the smaller crossing angle of $35^\circ$ a higher SNR of $20$ was assumed. In this case a SH-representation of $L=8$ is nearly 
at its limits to discriminate between both directions, while the discrete version can still well distinguish. Further, one can observe that for both methods the FODs along 
the horizontal Cartesian axis are sharper than along the skew axis. But this effect is more prominent for the discretized version than for the spherical harmonic representation. 

In Figure \ref{fig:results} we show quantitative results. The crossing was simulated for varying crossing angles between $30^\circ$ and $90^\circ$. Additionally we varied 
the absolute pose $\alpha$ of the crossing. For $\alpha=0^\circ$ the horizontal bundle is aligned with underlying x-axis of the Cartesian grid. With growing $\alpha$ the whole configuration 
is rotated clockwise. The crossing was simulated at a SNR of $50$, which is a realistic scenario. As a baseline experiment we show results of the so called Constrained Spherical
Deconvolution (CSD) approach \cite{Tournier2007}, where an additional positivity constraint is used to obtain more stable results. 
Apart from the unregularized CSD approach the results obviously depend on the absolute angle of the configuration. The discrete approach is able to cope quite good with 
small crossing angles, but shows independent of the crossing angle less precision than the SH-based approach. In particular, for crossing angles above $45^\circ$ the 
SH-based approach solves the task nearly without any error. 

In Figure \ref{fig:scatter} the results are further investigated by scatter plots in the $\phi,\theta$-plane. With the same parameters as above a crossing 
of $50^\circ$ was simulated and reconstructed by the three methods. Each detected direction is indicated by a small dot at its corresponding angle $\phi,\theta$. 
The crossing was simulated twice, for an absolute angle of $\alpha=0$, i.e. one direction is along the x-axis, and secondly for an absolute angle of $\alpha = 15^\circ$.
For $\alpha=0$ the SH-based approach is able to resolve the direction along the x-axis ($\phi=90^\circ$) perfectly, while the other direction ($\phi = 40^\circ$) is a bit
more blurry. On the other hand, the discrete approach has severe problems with the $\phi=90^\circ$ direction, which explains the lack of precision, i.e. apart from the 
true direction there are some additional local maxima that produce false positives. The main reason  is the interplay of the $64$ gradient directions and the 
$512$ discrete directions of the FOD. The effect is reduced by an increase of measurement directions. But also the SH-based approach has problems when the 
number of measurement directions is too low. They are revealed for an absolute angle of $\alpha=15^\circ$. Besides the uncertainty caused by the measurement noise
one can observe a systematic bias. For example, for the direction along $\phi=55^\circ$, the center of the distribution is shifted in $\theta$ by 
approximately $5^\circ$. Also for the other direction the distribution is a bit squeezed. We also found that the main reason is the low number of measurement 
directions. For example, for $128$ gradient directions the estimated directions are unbiased. Another way to reduce the effect is to decrease the 
expansion cutoff of the spherical harmonic representation. For $L=6$ and $64$ gradient directions the estimates do not show a bias. 
To conclude the differences: both methods have problems when the number of measurements become too low. While the discrete approach shows scattered, 
multimodal distributions, the SH-based approach shows a slight systematic bias but the distributions stay unimodal. 

In Figure \ref{fig:invivo_glyph} we show a real world example of the human brain. The setting of the measurement is nearly the same like in the simulations. A b-value of 
$1000 s/mm^2$ and $61$ gradient directions were used with an isotropic resolution of $2mm$. 
A kernel of the form $S_{\mv n_\text{fib}}(\mv n) = e^{-bD(\mv n \cdot \mv n_\text{fib} )^2}$ with $bD = 1$ was used as a model for deconvolution. 
Overall, both methods work very similar. Figure \ref{fig:invivo_glyph} shows a coronal section in glyph representation. One can observe that the SH-based 
approach produces a bit more negative values, which are indicated in black. The green rectangle highlights a region where the differences are largest.
The red rectangle shows a regions where the discrete approach shows a direction which does not appear for the SH-based approach. Whether the direction 
is true or not is difficult to say, but the fact that it precisely points along the x-axis makes it dubious. 

\subsubsection{Memory Consumption and Running Time}
The memory consumption of the SH-based and discrete approach is easy to compare. We want to consider the above real world experiment as an example. 
The whole volume has a size of $96 \times 96 \times 60$. As the FODs are symmetric we need to store only the even expansion coefficients of 
the spherical harmonic representation. Additionally, the FODs are real, hence we have to store $j+1$ numbers per expansion field, instead of $2j+1$ 
like for a general complex field due to the symmetry $\conj{Y^j_m} = (-1)^m Y^j_{-m}$. Thus, each voxels consumes $(L+2)^2/4$ complex numbers, resulting in $50\cdot 8$ bytes for $L=8$ and double precision. 
Overall, one volume needs $96^2 \cdot 60 \cdot 50 \cdot 8 = 220$ MB in SH-representation. On the other hand, in discrete representation with $512$ direction 
needs $96^2 \cdot 60 \cdot 512 \cdot 8 =  2264$ MB, which is $10$-times more compared to the SH-representation. Recall, that the conjugate gradient algorithm needs 
four instances of the volume in memory. 

The running time is more difficult to compare. On the one hand, we have to compute the application of $\Tk{0}$ which basically consists of finite differences,
on the other hand, the convolution operator $\mv H$ has to be implemented. 
Let us consider the computation of $\Tk{0}$ first. In case of the discrete approach one has to compute for each of the $512$ components 
six second order finite differences, which are linearly combined with weights depending on the corresponding direction. 
For the SH-based approach one also have to compute six finite differences, but the linear combinations of them to obtain the final values are more 
more expensive. In particular, we implemented separated functions for the operators $\mv \ME^2_{j,j+2}$, $\mv \ME^2_{j,j-2}$ and $\mv \ME^2_{j,j}$,
thus several values are computed repeatedly. In practice we found that one application of $\Tk{0}$ with $256$ discrete directions is comparable 
to an SH-based application of $\Tk{0}$ with $L=10$. 

The computation of the operator $\mv H$ is negligible in SH-representation, while it is the bottleneck for 
the discrete approach. Here, the running time heavily depends on the algorithm used for the matrix multiplication. In fact, the execution 
times can differ about a factor of $10$. While the highly optimized BLAS matrix-multiplication shipped with MATLAB is quite competitive, a standard non-optimized 
version can slow down the running time dramatically. To give an example, to reconstruct the above volume with $100$ CG iteration with $512$ directions 
takes on a \emph{Intel Xeon X7560} @ 2.27GHz about 20 minutes with a highly optimized multiplication. On the other hand, with a standard BLAS implementation it needs above an hour. 
For comparison, our implementation for $L=8$ takes about $10$ minutes on the same machine.

\subsection{A Spherical Hough Transform}

\begin{figure}
\centering
\includegraphics[width=16cm]{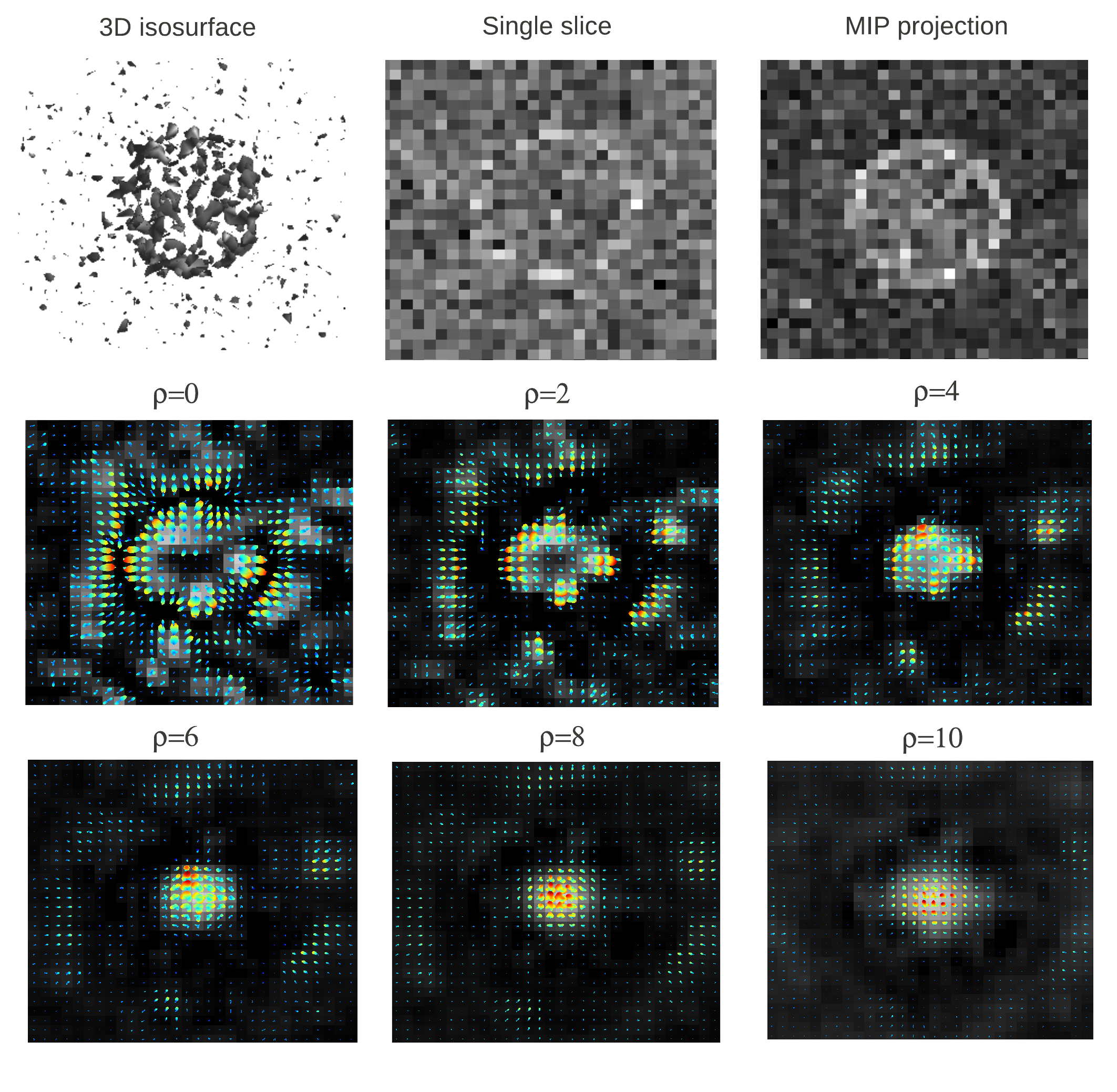}
\caption{A toy example for the spherical Hough transform. A distorted shell was rendered into a volume. The top row shows the 
initial data. In the two bottom rows the initial gradient orientation distribution is integrated with the kernel $\Tk{0} + 0.1\cdot\Tk{0} ^2$, with 
a step width of $\Delta\rho = 0.1$. The orientation distribution is expanded up to $L=4$. The gray value background shows the $j=0$ component.} \label{fig:houghtoy}
\end{figure}

\begin{figure}
\centering
\includegraphics[width=16cm]{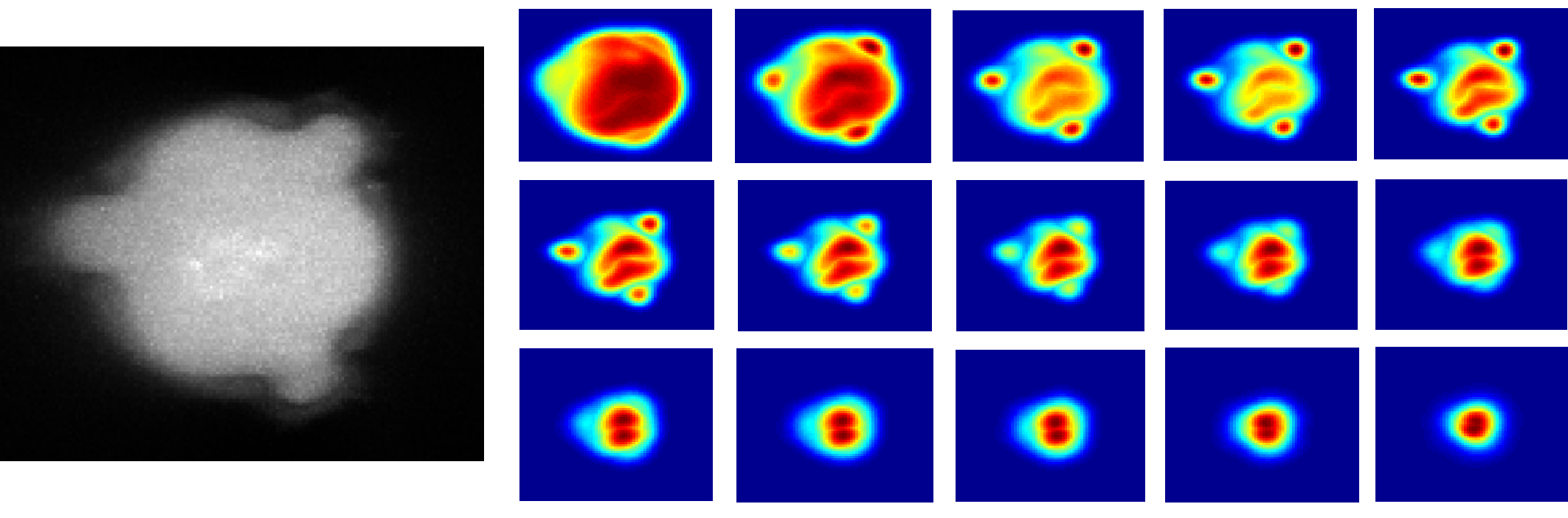}
\caption{The evolution of the Hough transform for a real world example of confocal laser scan of airborne pollen grain. All images are maximum intensity projections of 
the 3D volume image. The Hough transform was applied with the same setting like in the toy example.} \label{fig:houghreal}
\end{figure}

The detection of spherical objects is an everyday issue in biological image processing. The Hough transform \cite{hough} is often the method of choice. 
In this section we want to show how the contour completion kernel  \cite{duits:leftinv} can be used to implement a variant of this approach. Suppose, we have a volumetric image 
of a solid spherical object, such that the gradient field $\mv g(\mv r)$ at the surface of the object points away from the object center $\mv r_0$.
Let $m(\mv r) = |\mv g(\mv r)|$ be the gradient magnitude and $\mv v = \mv g/m$ the gradient direction. 
 If the object has radius $\rho$, then we know that for all $\mv n\in S_2$ approximately $\mv v(\mv r_0 + \rho\mv n) = \mv n$ holds. 
This fact can be easily used to get a 'evidence' or 'voting' map for the center $\mv r_0$ of the object. Therefore, let $\delta_{\mv n}(\mv m)$ 
an indicator function on the two-sphere, not necessarily a delta-function but a bit blurred. Then, we count for each putative 
object center $\mv r_0$ how often $\mv v(\mv r_0 + \rho\mv n) = \mv n$  is approximately fulfilled for each possible direction $\mv n \in S_2$:
\begin{eqnarray*}
h(\mv r_0, \rho) &=&   \int_{S_2}\ m(\mv r_0 + \rho \mv n) \ \delta_{\mv n}(\mv v(\mv r_0+\rho \mv n))\ d\mv n,
\end{eqnarray*}
where each contribution is weighted by the gradient magnitude of its originating voxel. The local maxima of the map $h$ give evidence for 
the presence of a spherical object with center $\mv r_0$ and radius $\rho$. 
Now, the integrand of the above equation can be expressed in terms of the horizontal translation operator:
\begin{eqnarray*}
H(\mv r,\mv n, \mv \rho) &=&   m(\mv r + \rho \mv n) \ \delta_{\mv n}(\mv v(\mv r+\rho \mv n)) \\
&=&  e^{-\rho\mv n \cdot \nabla} \left(  m(\mv r) \ \delta_{\mv n}(\mv v(\mv r)) \right) = e^{-\rho\Tk{0}} \left(  m(\mv r) \ \delta_{\mv n}(\mv v(\mv r)) \right)
\end{eqnarray*}
That is, we have a quite simple algorithm to get the Hough voting map: we initialize with $H(\mv r,\mv n,0) :=  m(\mv r) \delta_{\mv n}(\mv v(\mv r))$
and successively Euler-integrate 
\[
H(\mv r,\mv n,\rho^{(n+1)}) = H(\mv r,\mv n,\rho^{(n)}) + \Delta\rho\ \mathcal{A} H(\mv r,\mv n,\rho^{(n)}),
\]
where $\mathcal A = -\Tk{0}$. To get a more stable response we added, as discussed above, a slight amount of diffusion $\mathcal{A} = -\Tk{0} + 0.1\cdot\Tk{0} ^2$. 
Note that the proposed Hough transform only works for solid objects with surface gradients pointing away from the center. However, it is easy to switch
to inward pointing gradients by using  $\mathcal{A} = \Tk{0} + 0.1\cdot\Tk{0} ^2$ instead.

In Figure \ref{fig:houghtoy} we show a toy example for a shell, that is, the voxel on the surface of the sphere are set to one. 
Hence, both gradients are present: inward pointing gradients at the outer border and outward pointing gradients at the inner border. 
We decided to let the outward pointing gradients to be translated towards the center, i.e.  $\mathcal{A} = \Tk{0} + 0.1\cdot\Tk{0} ^2$. 
The shell in Figure \ref{fig:houghtoy}
was rendered with a radius of $7$ in a $32^3$-voxel cube. Every second voxel on the surface was deleted randomly and Gaussian random noise with standard deviation 
of $0.3$ was added. To get an impression Figure \ref{fig:houghtoy} shows an isosurface at gray value level $0.6$, a central slice and a maximum intensity projection (MIP)
of the toy example. The orientation field $H$ is expanded up to $L=4$ and integrated with a $\Delta\rho = 0.1$.
The initial gradient field $\mv g$ was computed on a Gaussian smoothed image of width $\sigma=1$. 
The lower two rows in Figure \ref{fig:houghtoy} show the evolution of orientation field $H$ for $\rho=0,2,4,6,8,10$ in glyph representation. The underlying gray value
image is the final voting map $h$, i.e. the $\ell=0$ component of the orientation field $H$. 
One can see how the gradients at the outer border are shifted outwards and gradients of the inner border are translated towards the center. 
Approximately for $\rho=7$ the voting maps looks as desired. 


\section{Conclusion}
This article worked out the formulation left-invariant convection/diffusion equations on $SE(3)$ in terms of the irreducible 
representations of $SO(3)$. From a computational viewpoint the main advantage of the proposed formulation is the low memory consumption
in comparison to a angular discretization of the two-sphere or rotation group, respectively. With a low number of basis functions the 
functions are well described without hurting the rotation covariance. The DWI example showed that even with 10 times less memory consumption 
the results are still comparable to the angular discrete approach. In terms of accuracy both approaches have their own problems. 

Applications to the full group $SE(3)$ remain subject to future work. For example, the detection of helical structures in cryo electron micro-graphs
\cite{ma:2012} might be a good playground. We further plan to provide the elementary operators as a open source toolbox to give the scientific community
a chance to try the proposed framework with a small amount of effort.

\section{Appendix}

\subsection{Spherical Harmonics} \label{ap:sh}
We always use Racah-normalized spherical harmonics such that $\mv Y^\ell(\mv r)^\T\mv Y^\ell(\mv r) = 1$, or
$\mv Y^\ell(\mv r)^\T\mv Y^\ell(\mv r') = P_\ell(\cos(\mv r,\mv r'))$, where the $P_\ell$ are the
Legendre polynomials:
\[
P_\ell(t) = \frac{1}{2^\ell \ell!} \partial_t^\ell (t^2-1)^\ell.
\]
 In terms of the associated Legendre polynomials
the components $Y^\ell_m$ of the spherical harmonics are written as
\[
Y^\ell_m (\phi,\theta) = \sqrt{\frac{(l-m)!}{(l+m)!}} P_\ell^m(\cos(\theta)) e^{\im m\phi}
\]
Mostly we write $\mv r \in S^2$ instead of $(\phi,\theta)$. The Racah-normalized
solid harmonics can be written as
\[
R^\ell_m(\mv r) = \sqrt{(\ell+m)!(\ell-m)!} \sum_{i,j,k} \frac{\delta_{i+j+k,\ell}\delta_{i-j,m}}{i!j!k! 2^i 2^j}
(x-\im y)^j (-x-\im y)^i z^k,
\]
where $\mv r = (x,y,z)$. They are related to spherical
harmonics by $R^\ell_m(\mv r) /r^\ell = Y^\ell_m(\mv r)$

\begin{eqnarray}
\langle J 0 | \ell_1 0,\ell_2 0 \rangle  Y^J_M = \sum_{m_1,m_2} \langle J M  | \ell_1 m_1,\ell_2 m_2 \rangle  Y^{\ell_1}_{m_1}  Y^{\ell_2}_{m_2}  \label{eq:shprod1} \\
\langle \ell_2 0 | \ell_1 0,J 0 \rangle  Y^J_M = \sum_{m_1,m_2} \langle \ell_2 m_2  | \ell_1 m_1, J M \rangle  \conj{Y^{\ell_1}_{m_1}}  Y^{\ell_2}_{m_2} \label{eq:shprod2} \\
Y^{\ell_1}_{m_1}  Y^{\ell_2}_{m_2} = \sum_{J,M} \langle J M | \ell_1 m_1,\ell_2 m_2 \rangle \langle J 0 | \ell_1 0,\ell_2 0 \rangle Y^{J}_{M} \label{eq:shprod3}\\
\conj{Y^{\ell_1}_{m_1}}  Y^{\ell_2}_{m_2} = \sum_{J,M} \frac{2J+1}{2\ell_1+1} \langle \ell_1 m_1 | \ell_2 m_2,J M \rangle \langle \ell_1 0 | \ell_2 0,J 0 \rangle Y^{J}_{M}\label{eq:shprod4} 
\end{eqnarray}

\subsection{Clebsch Gordan Coeffcients} \label{ap:cg}
The Clebsch Gprdan coefficients of $SO(3)$ fulfill several orthogonality relations:
\begin{eqnarray}
\sum_{j,m}  \langle j m | j_1 m_1, j_2 m_2 \rangle \langle j m | j_1 m'_1, j_2 m'_2 \rangle
&=& \delta_{m_1,m_1'} \delta_{m_2,m_2'}  \label{cg:orth1} \\
\sum_{j,m} \frac{2j+1}{2j_1+1} \langle j_1 m_1 | j m, j_2 m_2 \rangle \langle j_1 m_1' | j m, j_2 m_2' \rangle
&=& \delta_{m_1,m_1'} \delta_{m_2,m_2'}  \label{cg:orth1b} \\
\sum_{{m = m_1 + m_2} } \langle j m | j_1 m_1, j_2 m_2 \rangle \langle j' m' | j_1 m_1, j_2 m_2 \rangle
&=& \delta_{j,j'} \delta_{m,m'}\label{cg:orth2}  \\
\sum_{m_1,m } \langle j m | j_1 m_1, j_2 m_2 \rangle \langle j m | j_1 m_1, j'_2 m'_2 \rangle
&=& \frac{2j+1}{2j'_2+1} \delta_{j_2,j'_2} \delta_{m_2,m_2'}\label{cg:orth3} 
\end{eqnarray}
For particular combinations there are simple formulas:
\begin{eqnarray}
\langle \ell m | (\ell-\lambda)(m-\mu), \lambda\mu \rangle =
\left(\begin{array}{c} \ell+m \\ \lambda + \mu \end{array} \right)^{1/2}
\left(\begin{array}{c} \ell-m \\ \lambda - \mu \end{array}\right)^{1/2}
\left(\begin{array}{c} 2\ell \\ 2\lambda  \end{array}\right)^{-1/2}
\end{eqnarray}
\begin{eqnarray}
\begin{split}
\langle \ell m | (\ell+\lambda)(m-\mu), \lambda\mu \rangle& = (-1)^{\lambda+\mu}
\left(\begin{array}{c} \ell+\lambda -m + \mu \\ \lambda + \mu \end{array} \right)^{1/2} \\
&\left(\begin{array}{c} \ell+\lambda+m-\mu \\ \lambda - \mu \end{array}\right)^{1/2}
\left(\begin{array}{c} 2\ell + 2\lambda + 1 \\ 2\lambda  \end{array}\right)^{-1/2}
\end{split}
\end{eqnarray}
There are several symmetry relations
\begin{eqnarray}
 \langle j m | j_1 m_1, j_2 m_2 \rangle &=& \langle j_1 m_1, j_2 m_2| j m \rangle \\
 \langle j m | j_1 m_1, j_2 m_2 \rangle &=& (-1)^{j+j_1+j_2} \langle j m | j_2 m_2, j_1 m_1 \rangle \\
 \langle j m | j_1 m_1, j_2 m_2 \rangle &=& (-1)^{j+j_1+j_2} \langle j (-m) | j_1 (-m_1), j_2 (-m_2) \rangle \\
 \langle j m | j_1 m_1, j_2 m_2 \rangle &=& \sqrt{\frac{2j+1}{2j_2+1} }(-1)^{j_1+m_1} \langle j_2 m_2 | j m, j_1 (-m_1) \rangle ,
\end{eqnarray}
and associativity relations:
\begin{eqnarray}
&\langle J,M|j_1+j_2,m_1+m_2,j_3,m_3\rangle  \langle j_1+j_2,m_1+m_2|j_1,m_1,j_2,m_2\rangle = \nonumber\\
&\quad \quad \quad \langle J,M|j_1+j_3,m_1+m_3,j_2,m_2\rangle  \langle j_1+j_3,m_1+m_3|j_1,m_1,j_3,m_3 \rangle
\end{eqnarray}
where $J=j_1+j_2+j_3$ and $M=m_1+m_2+m_3$. And for $j_3 > j_1+j_2$ we have another one:
\begin{eqnarray}
&\langle j_3-j_1-j_2,m_1+m_2+m_3|J-j_1,m_1+m_3,j_2,m_2\rangle  \langle j_3-j_1,m_1+m_3|j_1,m_1,J,m_3\rangle = \nonumber\\
&\langle j_3-j_1-j_2,m_1+m_2+m_3|J-j_2,m_2+m_3,j_1,m_1\rangle  \langle j_3-j_2,m_2+m_3|j_2,m_2,J,m_3\rangle = \nonumber\\
\end{eqnarray}

\subsection{Wigner D-Matrix} \label{ap:wigner}
The irreducible representation of $SO(3)$ are called Wigner D-matrices $\mRep{g}^\ell$ and are indiced by an integer $\ell=0,\hdots,\infty$. The $\ell$th
order representation works on a $\Complex^{2\ell+1}$ dimensional vector space. We denote the components of $\mRep{g}^\ell$ by $D^\ell_{mn}(g)$. In Euler angles
in ZYZ-convention we have
\begin{eqnarray}
D^\ell_{mn}(\gamma,\beta,\alpha) = e^{-\im m\gamma} d^\ell_{mn}(\beta) e^{\im n \alpha},
\end{eqnarray}
where $d^\ell_{mn}(\beta)$ is the 'small' Wigner d-matrix, which is real-valued and explicitly written as
\begin{eqnarray}
\begin{array}{lcl}
d^\ell_{mn}(\beta) = [(\ell+m)!(\ell-m)!(\ell+n)!(\ell-n)!]^{1/2}
\sum_s \frac{(-1)^{m-n+s}}{(\ell+n-s)!s!(m-n+s)!(\ell-m-s)!} \\
\quad \quad \times \left(\cos\frac{\beta}{2}\right)^{2\ell+n-m-2s}\left(\sin\frac{\beta}{2}\right)^{m-n+2s}.
\end{array}
\end{eqnarray}
The representations of different order are connected via the Clebsch Gordan coefficients by:
\begin{eqnarray}
D^\ell_{mn} = \sum_{m_1+m_2 = m \atop n_1+n_2 = n} D^{\ell_1}_{m_1 n_1}D^{\ell_2}_{m_2 n_2}
 \langle l m | l_1 m_1, l_2 m_2 \rangle
 \langle l n | l_1 n_1, l_2 n_2 \rangle \label{dw:decomp1}
\end{eqnarray}
and
\begin{eqnarray}
D^{\ell_1}_{m_1 n_1}D^{\ell_2}_{m_2 n_2} = \sum_{l,m,n} D^{\ell}_{m n}
 \langle l m | l_1 m_1, l_2 m_2 \rangle
 \langle l n | l_1 n_1, l_2 n_2 \rangle \label{dw:decomp2}
\end{eqnarray}
Another important equality is 
\begin{eqnarray}
\int_{SO(3)} dg\ \conj{D^\ell_{k' k}} \conj{D^{j'}_{n' m'}}  D^{j}_{n m} = \frac{8 \pi^2}{2 j + 1}  \langle j n | j' n', \ell k' \rangle \langle j m | j' m' ,\ell k \rangle \label{eq:tripleproducts}
\end{eqnarray}
%
%

\subsection{Mixed Quadratic Terms} \label{ap:mixedquad}
The action of the terms $\Tk{\mp} \Jpm$ and $\Jpm \Tk{\mp} $ can be computed directly from the equations \eqref{eq:DwigQM1}, \eqref{eq:DwigQM2}, \eqref{eq:transfield2} 
to
\begin{eqnarray*}
(\Tk{\mp} \Jpm \mv f)^j_{nm} &=&  -\im
  \sum_{j'=j-1,j,j+1\atop q  = -1,0,1} \sum_{n= n'+q } \sqrt{j'(j'+1)/2-(m\pm 1)(m\pm 2)/2} \frac{2j'+1}{2 j + 1}  \\ 
&&\langle j n | j' n', 1 q \rangle \langle j m | j' (m\pm 1) ,1 (\mp 1) \rangle\ \gmv \partial^1_{q} f^{j'}_{n'm\pm 2}
\end{eqnarray*}
and
\begin{eqnarray*}
(\Jpm \Tk{\mp} \mv f)^j_{nm} &=&  \im \sqrt{j(j+1)/2-m(m\pm 1)/2}  \\
&&\sum_{j'=j-1,j,j+1\atop q  = -1,0,1} \sum_{n= n'+q \atop m = m'  \mp 2}  \frac{2j'+1}{2 j + 1}  \langle j n | j' n', 1 q \rangle \langle j (m\pm 1) | j' (m\pm 2) ,1 (\mp 1) \rangle\ \gmv \partial^1_{q} f^{j'}_{n'm \pm 2}
\end{eqnarray*}

\subsection{Proof of equation \eqref{eq:shquadrel}} \label{ap:proofshquad}
From equation \eqref{eq:shprod4}  we directly know the corresponding equality for the spherical derivatvies:
\[
\conj{\gmv \partial ^{\ell_1}_{m_1}} \gmv \partial^{\ell_2}_{m_2} = \sum_{J,M} \frac{2J+1}{2\ell_1+1} \langle \ell_1 m_1 | \ell_2 m_2,J M \rangle \langle \ell_1 0 | \ell_2 0,J 0 \rangle\gmv \partial^{J}_{M}
\]
Setting $\ell_1 = \ell_2 = 1$ we know that the sum over $J$ takes only values for $J=0$ and $J=1$. By rotating the frame of reference 
$\gmv \partial^\ell \mapsto  \mv D^\ell(g)^\T \gmv \partial^\ell$ we get 
\[
(\conj{\mv D^1(g)^\T \gmv \partial ^1})_{m_1} (\mv D^1(g)^\T \gmv \partial^{1})_{m_2} = \sum_{J,M} \frac{2J+1}{3} \langle 1 m_1 | 1 m_2,J M \rangle \langle 1 0 | 1 0,J 0 \rangle (\mv D^J(g)^\T\gmv\partial^{J})_{M}
\]
with $\Tk{m} =  (\mv D^1(g)^\T \gmv \partial^{1})_{m}$ and evaulting the Clebsch Gordan coefficients we end up with 
\begin{equation*} 
\conj{\Tk{m_1}} \Tk{m_2} =  \frac{\Delta}{3}  -  \frac{\sqrt{10}}{3} \sum_{M=-2}^2  \langle 1 m_1 | 1 m_2,2 M \rangle  \ (\mv D^2(g)^\T \gmv \partial^2)_M
\end{equation*}
which was to show.

\bibliographystyle{plain}

\end{document}